\documentclass[11pt]{article}
\newcommand{\version}{September 18, 2009}

\usepackage{amsthm,amsfonts,amsmath,amssymb}
\usepackage{graphicx}
\usepackage{pdfsync}
\usepackage{appendix}

\swapnumbers
\theoremstyle{plain}

\newtheorem{thm}{THEOREM}[section]
\newtheorem{lm}[thm]{LEMMA}

\newtheorem{prop}[thm]{PROPOSITION}
\newtheorem{remark}[thm]{REMARK}

\theoremstyle{definition}
\newtheorem{defi}[thm]{DEFINITION}
\theoremstyle{remark}
\newcommand{\upchi}{\raise1pt\hbox{$\chi$}}
\newcommand{\R}{{\mathord{\mathbb R}}}

\newcommand{\lanbox}{\hfill \hbox{$\, \vrule height 0.25cm width 0.25cm depth 0.01cm\,$}}

\def\d{{\rm d}}

\numberwithin{equation}{section}  

\begin{document}

\markboth{\scriptsize{CGL \version}}{\scriptsize{CGL \version}}

\title{\bf{ On the Markov sequence problem for Jacobi polynomials}}
\author{\vspace{5pt} Eric A. Carlen$^{1,3}$,  Jeffrey S. Geronimo$^{2,4}$, and Michael Loss$^{2,5}$ \\
\vspace{5pt}\small{$1.$ Department of Mathematics, Hill Center, Rutgers University}\\[-6pt]
\small{
110 Frelinghuysen Road
Piscataway NJ 08854 USA}\\
\vspace{5pt}\small{$2.$ School of Mathematics, Georgia Institute of
Technology,} \\[-6pt]
\small{Atlanta, GA 30332 USA}\\
 }
\date{\version}
\maketitle \footnotetext [3]{Work of Eric Carlen is partially supported by U.S.
National Science Foundation
grant DMS-0901632}
\maketitle \footnotetext [4]{Work of Jeffery Geronimo is partially supported by U.S.
National Science Foundation
grant DMS-0500641.}  
\maketitle \footnotetext [5]{Work of Michael Loss is partially supported by U.S.
National Science Foundation
grant DMS0-901304 \\
\copyright\, 2008 by the authors. This paper may be reproduced, in
its entirety, for non-commercial purposes.}

\begin{abstract}  We give a simple and entirely elementary
proof of Gasper's theorem on the Markov sequence problem for Jacobi polynomials.
It is based on the spectral analysis of an operator that arises in the study of a probabilistic model of colliding molecules introduced by Marc Kac, and the methods developed here yield new estimates relevant to
the collision model.

\end{abstract}

 \centerline{Mathematics Subject Classification Numbers: 31B10, 33C45, 37A40}

\section{Introduction} \label{MS}

\subsection{The Markov Sequence Problem and the Theorems of Bochner and Gasper}

Let $(X, {\mathcal  S},\mu)$ be a probability space. A {\it Markov operator} $T$  on $L^2(\mu)$  is a linear operator 
 that preserves positivity; i.e., $f\ge 0 \Rightarrow Tf \ge 0$, and preserves the constants; i.e., $T1 =1$. 
If $T$ is self adjoint, it follows by duality and the Riesz-Thorin interpolation theorem that T is a
contraction  on $L^p(\mu)$ for all $1\le p \le
\infty$,)
Consequently, the spectrum of $T$ lies in the interval $[-1,1]$. 

The next definitions, which are less standard, are taken from \cite{BH} and \cite{BM}:  A {\it unit orthonormal basis} for 
$L^2(\mu)$ is an orthonormal basis   $\{f_n\}_{n\ge 0}$  such that $f_0 =1$.  Though we discuss a broader class of examples in Section 5 and in the Appendix, In the main examples  here, 
 $X= \R$, or some subset of $\R$, and 
$\{f_n\}_{n\ge 0}$
is the sequence of orthonormal polynomials for $\mu$. In any case, we shall always suppose that $X$
is a locally compact Hausdorff space, and that $\mu$ is a Borel  measure. 

Given a  unit orthonormal basis $\{f_n\}_{n\ge 0}$, the set of  {\it Markov sequences}  ${\mathcal  M}$
for this basis is the set of all sequences  $\{\lambda_n\}_{n\ge 0}$ such that there exists a self adjoint Markov operator $K$ with
$$Kf_j = \lambda_j f_j\qquad{\rm for\ all}\quad j\ge 0\ .$$
Notice that necessarily $\lambda_0 =1$ and $\lambda_n\in [-1,1]$ for all $n$. Also, since a convex combination of self adjoint Markov operators is self adjoint and Markov, ${\mathcal  M}$ is convex, so that
 ${\mathcal  M}$ may be described by specifying its extreme points. 

The  {\it Markov sequence problem} is to determine, for a given unit orthonormal basis, the set ${\mathcal  M}$. Naturally, it is sufficient
to find the extreme points.

The  Markov sequence problem seems to have been first considered by  Bochner \cite{B}, and the first result, for ultraspherical polynomials, is his as well. 

We recall that the for each $\gamma>-1/2$, the {\em ultraspherical polynomials} 
$\{p^{(\gamma)}_n\}_{n\ge 0}$ are the {\em orthonormal polynomials},  for the measure
$\mu^{(\gamma)}$ 
\begin{equation}\label{mualpha}
{\rm d}\mu^{(\gamma)}(t) =  c_\gamma (1-t^2)^{\gamma-1/2}{\rm d}t \quad{\rm where}\quad
c_\gamma = \frac{1}{\sqrt{\pi}}\frac{\Gamma(\gamma+1)}{\Gamma(\gamma+1/2)}
\end{equation}
is the normalization constant that makes $\mu^{(\gamma)}$ a probability measure. 
The normalization as unit vectors in $L^2(\mu^{(\gamma)})$ is just one useful and frequently encountered normalization. Another that will be useful here is 
generally denoted with an upper--case $P$: 
The ultraspherical polynomials $P_n^{(\gamma)}$ are normalized so that  
$P_n^{(\gamma)}(1) =1$ i.e.
\begin{equation}\label{gegen2}
P_n^{(\gamma)}(x)  = \frac{p_n^{(\gamma)}(x)}{p_n^{(\gamma)}(1)} \ .
\end{equation}
  Throughout the paper, an upper--case $P$ denotes this normalization, while a lower case $p$ denotes the $L^2(\mu^{(\gamma)})$ normalization.

The ultraspherical polynomials are special cases in the wider family of Jacobi polynomials: Recall that the Jacobi polynomials $p_n^{(\alpha,\beta)}$ form an
orthonormal basis for
$L^2\left([-1,1], d\mu^{\alpha,\beta}\right)$
where 
\begin{equation}\label{mabdef}
\mu^{\alpha,\beta}({\rm d}x) = c_{\alpha,\beta}(1-x)^{\alpha}(1+x)^{\beta}{\rm d}x \ ,
\end{equation}
where $c_{\alpha, \beta}$ makes $\mu^{\alpha,\beta}$ a probability measure. 
In particular, the ultraspherical polynomials arise for the special case $\gamma = \alpha-1/2 = \beta-1/2$;
that is
\begin{equation}\label{convent}
p_n^{(\gamma)}(t) =   p_n^{(\gamma-1/2,\gamma-1/2)}(t) \ .
\end{equation}

Theorem 2 of \cite{B} may be phrased as follows:

\begin{thm}[\bf Bochner]\label{bochner}  For any  $\gamma> 0$,  the sequence $\{\lambda_n\}_{n\ge 0}$ is a Markov sequence for
$\{p^{(\gamma)}_n\}_{n\ge 0}$ if and only if there is a probability measure $\nu$ on $[-1,1]$
such that
\begin{equation}\label{ultra}
\lambda_n = \int_{-1}^1\frac{p_n^{(\gamma)}(t)}{p_n^{(\gamma)}(1)}{\rm d}\nu(t)\ .
\end{equation}
For each such Markov sequence $\{\lambda_n\}_{n\ge 0}$, the measure $\nu$ is unique. In other words, for each $t$, $\{p_n^{(\gamma)}(t)/p_n^{(\gamma)}(1)\}_{n\ge 0}$ is a Markov sequence
for $\{p^{(\gamma)}_n\}_{n\ge 0}$, and these are the extreme points of the set ${\mathcal  M}$ of all such Markov sequences.
\end{thm}

Since the ultraspherical polynomials are Jacobi polynomials with $\alpha = \beta$, 
it is natural to ask whether one can one extend Bochner's result  to a wider class of Jacobi polynomials with $\alpha \ne \beta$. 
This  question was answered by Gasper \cite{G1,G2}:

\begin{thm}[{\bf Gasper}]\label{gasp1}  For   $\alpha\ge \beta$ with  $\beta> -1/2$ or
$\alpha > \beta$ with  $\beta = -1/2$
the sequence $\{\lambda_n\}_{n\ge 0}$ is a  Markov sequence for
$\{p^{(\alpha,\beta)}_n\}_{n\ge 0}$, if and only if there is a probability measure $\nu$ on $[-1,1]$
such that
\begin{equation}\label{jac1}
\lambda_n = \int_{-1}^1\frac{p_n^{(\alpha,\beta)}(x)}{p_n^{(\alpha,\beta)}(1)}{\rm d}\nu(x)\ .
\end{equation}
For each such Markov sequence $\{\lambda_n\}_{n\ge 0}$, the measure $\nu$ is unique. In other words, for each $t$, $\{p_n^{(\alpha,\beta)}(t)/p_n^{(\alpha,\beta)}(1)\}_{n\ge 0}$ is a Markov sequence
for $\{p^{(\gamma)}_n\}_{n\ge 0}$, and these are the extreme points of the set ${\mathcal  M}$ of all such Markov sequences.
\end{thm}

While Bochner's proof of  Theorem~\ref{bochner} is not terribly complicated, Gasper's  proof of 
Theorem~\ref{gasp1} is far from elementary. Even though it has been simplified
by the work of others, particularly Koornwinder, it remains a {\em tour de force}: Koornwinder's proof still uses many deep results on special functions.  

In this paper we shall give entirely elementary and self-contained  proofs
of these theorems. Moreover, these proofs will allow us to obtain bounds on the sizes of the extremal eigenvalues. Before stating the new results more precisely, we recall the proof of Bochner's Theorem, as this will
clarify the matter of what was already well understood, and what was in need of clarification.

\subsection{Product formulas and the Markov sequence problem}

In this subsection we explain that the Markov sequence problem is easily solved for unit orthonormal
sequences that satisfy a {\em product formula}, as defined below.  Indeed, Bochner's original proof of his theorem went by this route, and was facilitated by the fact that the product formula he required had already been established long ago by Gegenbauer. Gasper, on the other hand, had more work to do since before his work, no general product formula for Jacobi polynomials was known.

\begin{defi}[Product Formula] A unit orthonormal sequence $\{f_n\}_{n\ge 0}$ satisifes a 
{\em product formula} in case there exists function $(x,y) \mapsto {\rm d}\mu_{x,y}(z)$
from $X\times X$ to the space of probability measures on $X$, 
and also some $x_0\in X$ such that for each $n\ge 0$,
\begin{equation}\label{product}
F_n(x)F_n(y) = \int_X F_n(z){\rm d}\mu_{x,y}(z)\ ,
\end{equation}
where ${\displaystyle F_n(x) = \frac{f_n(x)}{f_n(x_0)}}$.
\end{defi}

For example, in the case of the ultraspherical polynomials $\{p_n^{(\gamma)}\}_{n\ge 0}$, 
take $x_0= 1$, so that $F_n$ becomes
$P_n^{(\gamma)}$.  Then one has {\em Gegenbauer's identity} \cite{Geg}, which dates back to 1875:

\begin{thm}[\bf Gegenbauer's Identity]\label{gegid} For all $\gamma>1/2$, and all $n\ge 0$, and all
$a\in (-1,1)$, 
\begin{equation}\label{boc12}
P^{\gamma}_n(a)P^{\gamma}_n(t) = \int_{-1}^1  P^{\gamma}_n\left(at   + s \sqrt{1-a^2}\sqrt{1-t^2}\right){\rm d}\mu^{(\gamma-1/2)}(s)\ .
\end{equation} 
\end{thm}

To see this as a concrete  instance of the abstract  product formula (\ref{product}), let $\delta_u$ denote the Dirac mass at $u\in [-1,1]$, and define
$${\rm d}\mu_{a,t}(z) = \int_{-1}^1  \delta_{at   + s \sqrt{1-a^2}\sqrt{1-t^2}}(z){\rm d}\mu^{(\gamma-1/2)}(s)\ .$$
Then (\ref{boc12}) becomes ${\displaystyle 
P^{\gamma}_n(a)P^{\gamma}_n(t) = \int_{-1}^1  P^{\gamma}_n(z){\rm d}\mu_{a,t}(z)}$,
as in (\ref{product}).

The following theorem relates the Markov sequence problem to the problem of establishing a product formula. The theorem summarizes ideas that can be found, reading between the lines,  in Bochner's paper \cite{B} for the ultraspherical polynomials, and much more explicitly, and in general, in the paper \cite{BH} of Bakry and Huet.

\begin{thm}[\bf Markov Sequences and Product Formulae]\label{BHT}  Let $X$ be a closed interval in $\R$, and let $\mu$ be a regular Borel probability measure whose support is $X$. Let $\{f_n\}_{n\ge 0}$
be a unit orthonormal basis for 
$L^2(\mu)$ consisting of real valued functions.  Let $x_0$ be any fixed point in $X$. Then the following are equivalent:

\smallskip
\noindent{\em (1)} For each $x\in X$,  
$\{\lambda_n(x)\}_{n\ge0}$
is a Markov sequence for  $\{f_n\}_{n\ge 0}$ where
${\displaystyle \lambda_n(x) := \frac{f_n(x)}{f_n(x_0)}}$.

\smallskip
\noindent{\em (2)}  For each  $\{\lambda_n\}_{n\ge 0} \in {\mathcal  M}$, there exists a 
a  Borel probability measure
$\nu$ so that
\begin{equation}\label{ultra2}
\lambda_n = \int_X\frac{f_n(x)}{f_n(x_0)}{\rm d}\nu(x)\ .
\end{equation} 

\smallskip
\noindent{\em (3)} With $F_n(x) := f_n(x)/f_n(x_0)$, the $\{F_n\}_{n\ge0}$ satisfy the product formula
(\ref{product}) for some family ${\rm d}\mu_{x,y}(z)$ of probability measures on $X$.
\smallskip

Finally, if any (and hence all) of these conditions are satisfied, and if  $\{f_n\}_{n\ge 0}$ is a sequence of bounded continuous functions whose
finite linear combinations are dense in ${\mathcal  C}_b(X)$, then the probability measure $\nu$ in (\ref{ultra2}) is unique, so that ${\mathcal  M}$ is a simplex and the  $\{f_n(x)/f_n(x_0)\}_{n\geq 0}$ are   its extreme points.
\end{thm}

The equivalence of (1) and (2), as well as the statement concerning uniqueness of the measure $\nu$, is due to Bakry and Huet,  \cite{BH}, together with many other results on the Markov sequence problem.
The equivalence of (1) and (3) is implicit in Bocher's paper \cite{B}, though his argument is different from what follows below, and in particular, he makes no use of self-adjointness of certain operators associated to product formulae -- a crucial feature of our approach.  Thus, while we make no claim of originality for the results in Theorem~\ref{BHT}, 
 we provide a complete proof  for completeness and clarity.

\smallskip
\noindent{\bf Proof:} We first show that (1) implies (2): 
If $\{f_n(z)/f_n(x_0)\}_{n\geq 0}$ is in ${\mathcal  M}$, then 
by the spectral theorem, 
${\displaystyle K_z(x,y) = \sum_{n=0}^\infty \lambda_n(z) f_n(x) f_n (y)}$
is the kernel of a Markov operator $K_z$ on 
$L^2(\mu)$ with
$K_zf_n (x) =  \lambda_n(z)f_n(x)$. If $\nu$ is a Borel probability measure then
${\displaystyle 
K := \int_X K_z {\rm d} \nu(z)
}$
is a Markov operator with eigenvalues 
\begin{equation}
\int_X \lambda_n(z) {\rm d}\nu(z) = \int_X \frac{f_n(z)}{f_n(x_0)} {\rm d}\nu(z) \ .
\end{equation}

We next show that (2) implies (1): Suppose that $\{\lambda_n\}_{n \ge 0}$ is a Markov sequence. Again, by the spectral theorem ${\displaystyle K(x,y) := \sum_{n=0}^\infty \lambda_n f_n(x) f_n (y)}$
is the kernel of a Markov operator $K$ on 
$L^2(\mu)$ with
$Kf_n (x) =  \lambda_nf_n(x)$.
Since $K$ is Markov, $K(x,y) \ge 0$ for all $x$ and $y$, and for each $y$, $K(x,y){\rm d}\mu(x)$  is a probability measure. Taking $y=x_0$, define the probability measure
${\rm d} \nu =  K(x,x_0){\rm d}\mu$.
Then, for each $k$,
\begin{eqnarray}
\int_{X}f_k(x){\rm d} \nu(x) &=&  \int_{X}f_k(x)\left[\sum_{n=0}^\infty \lambda_n
f_n(x) f_n(x_0){\rm d}\mu\right]\nonumber\\
&=& \sum_{n=0}^\infty \left[ \int_{X}f_k(x) f_n(x){\rm d}\mu\right]\lambda_n
f_n(x_0)
\ =\ \lambda_kf_k(x_0)\ ,\nonumber\\
\end{eqnarray}
which proves (\ref{ultra2}).

We next show that (1) implies (3):  If $\{f_n(z)/f_n(x_0)\}_{n\geq 0}$ is in ${\mathcal  M}$, then 
for each $z$ in $X$,
$$
k_z(x,y) = \sum_{n=0}^\infty \frac{f_n(z) f_n(x) f_n (y)}{f_n(x_0)}
$$
is the eigenfunction expansion of the kernel of a self adjoint Markov operator. Evidently,
$(x,y,z)\mapsto k_z(x,y)$ is invariant under any permutation of $x$, $y$, and $z$. Thus, for each $m$,
\begin{equation}\label{prfo1}
\int_X f_m(z)k_z(x,y){\rm d}\mu(z) =  
\int_X k_x(y,z) f_m(z) {\rm d}\mu(z)  = \frac{f_m(x)}{f_m(x_0)}f_m(y)\ .
\end{equation}
To recognize this as a product formula, 
for each $x,y$, define a probability measure $\mu_{x,y}$ by
\begin{equation}\label{prfo1b}
{\rm d}\mu_{x,y}(z) = k_z(x,y){\rm d}\mu(z)\ .
\end{equation}
Then with 
$F_n(x) := f_n(x)/f_n(x_0)$,
(\ref{prfo1}) becomes  (\ref{product}).

We next show that (3) implies (1): Assuming (3), fix $y\in $ and define an operator $K_y$
on ${\cal C}_b(X)$ by 
$$K_y\varphi(x) = \int_X \varphi(z){\rm d}\mu_{x,y}(z)\ .$$
Since $K$ is a Markov operator, it has a bounded extension to $L^2(\mu)$.
The product formula says that for each $m$, $f_m$ is an eigenfunction of $K_y$
with eigenvalue $F_m(y)$. Any bounded operator with a complete orthonormal basis of eigenfunctions, each of whose eigenvalues is real, is necessarily self adjoint. Thus $K_y$ is a self adjoint Markov operator, and hence $\{F_n\}_{n\ge 0}$ belongs to ${\mathcal M}$.

Finally, It remains to show that
the measure $\nu$ is uniquely determined.
For this,  let $f$ be any continuous function on $X$. Let $\epsilon>0$ be given, and let
$g(x)= \sum_{n=0}^N\alpha_n f_n(x)$ be a finite linear combination of the
$\{f_n(x)/f_n(x_0)\}_{n\geq 0}$ such that $|f(x)  -  g(x)| \le \epsilon$
for all $x$. 

Let $\nu$ and $\widehat \nu$ be two Borel probabilty measures such that 
 (\ref{ultra2}) holds for some  $\{\lambda_n\}_{n\ge 0} \in {\mathcal  M}$. Then
 $$\int_X g(x){\rm d}\nu(x) = \sum_{n=0}^N\alpha_n\int_X f_n(x){\rm d}\nu(x) = 
 \sum_{n=0}^N\alpha_n\lambda_n f_n(x_0) = 
 \sum_{n=0}^N\alpha_n\int_X f_n(x){\rm d}\widehat  \nu(x) = 
 \int_X g(x){\rm d}\widehat  \nu(x)\ .$$
 Therefore, $|\int_X f(x){\rm d}  \nu(x) -  \int_X f(x){\rm d}\widehat  \nu(x)| \le 2\epsilon$.
Since $\epsilon$ is arbitrary, 
 $\int_X f(x){\rm d}  \nu(x) = \int_X f(x){\rm d}\widehat  \nu(x)$
 for all $f\in {\mathcal  C}_b(X)$. This of course means that  $\nu = \widehat \nu$.
 \lanbox

\begin{remark}\label{wiers}Notice that   by the Weierstrass Approximation Theorem, the conditions in the final part of 
Theorem~\ref{BHT} are automatically satisfied in  any  application to orthogonal polynomials on a compact interval. 
\end{remark}

It is interesting to note that once one has 
 the product formula (\ref{product}), one can use it to define a convolution: 
For any two finite, positive measures $\lambda, \nu$  on $X$,  define the {\em convolution} $\lambda \star \nu$
of $\lambda$ and $\nu$ by
\begin{equation} \label{convol}
\lambda \star \nu := \int {\rm d}\mu_{x,y} (z)  {\rm d}\lambda(x)  {\rm d}\nu(y) \ ;
\end{equation}
this too is a finite positive measure.  

Note that the ``Fourier'' coefficients of $\lambda$ and $\nu$, 
given by
${\displaystyle 
\int F_n  {\rm d}\lambda \ {\rm and} \ \int F_n  {\rm d}\nu}$
satisfiy
$$
\int F_n  {\rm d}\lambda \int F_n  {\rm d}\nu = \int F_n {\rm d}(\lambda \star \nu)  \ ,
$$
so that the usual relation between Fourier coefficients and convolutions holds.
We now return to the matter of proving Bochner's Theorem.

\if false
For more details, the reader may consult \cite{G1}, \cite{KooSchwa} and \cite{BH}, especially Section 2.4.
In Gasper's original development \cite{G1} of these ideas, the starting point for the definition of such a convolution
was to find a unit orthonormal bases $\{f_n\}_{n\ge 0}$ satisfying
\begin{equation}\label{prfo8}
\sum_{n=0}^\infty \frac{f_n(z) f_n(x) f_n (y)}{f_n(x_0)} \ge 0\ .
\end{equation}
Whenever this positivity condition is satisfied,  the sum on the left  is at least formally the eigenfunction expansion for the kernel of a Markov operator, and  after dealing with convergence issues,
 one can define the convolution operation as above. 
Convolutions are often associated with harmonic analysis on groups, but as
 indicated, the positivity condition (\ref{prfo8})allows one to define
 convolutions beyond this context. In Section 2.3  of \cite{BH}
 equation (\ref{prfo8}) is refered to as the hypergroup property.  A general discussion of hypergroups is beyond the
 scope of this paper and the interested reader may consult the book of Bloom
 and Heyer \cite{BlHe}
 \fi

\subsection{Gegenbauer's Identity and a proof of Bochner's Theorem}

\smallskip

\noindent{\bf Proof of Bochner's Theorem:} 
Since we have a product formula for the ultraspherical polynomials, namely Gegenbauer's identity (\ref{boc12}), condition (3) of Theorem~\ref{BHT} is satisfied, and Bochner's Theorem follows immediately.
\lanbox

\smallskip

This proof is simple, but hardly complete: one needs the product formula. In the case of the ultraspherical polynomials, this was ready at hand since 1875.  For the general case of the Jacobi polynomials, no product formula was available when Gasper began his work.   His strategy was to show that for
$\{f_n\}_{n\ge 0}$ being a sequence of Jacobi polynomials, one has the positivity result
\begin{equation}\label{prfo8}
\sum_{n=0}^\infty \frac{f_n(z) f_n(x) f_n (y)}{f_n(x_0)} \ge 0\ ,
\end{equation}
holding pointwise almost everywhere. Then, this sum defines a positive kernel, which can be
used, as in the proof of Theorem~\ref{BHT}, to prove a product formula. However, this direct proof of pointwise positivity is far from simple.

In this paper we present a truly simple approach to the product formula for Jacobi polynomials. First, however, we shall illustrate this approach by providing a simple, self-contained proof of Gegenbauer's  identity  (\ref{boc12}).

For our purposes, it is most  helpful to consider  (\ref{boc12})
as an eigenvalue identity.

\begin{defi}[\bf The Correlation Operators] \label{cor}
For each $\gamma>0$,  and $a\in (-1,1)$, define an operator $K_a$ on $L^2(\mu^{(\gamma)})$ by
\begin{equation}\label{Kdef}
K_af(t) = \int_{-1}^1  f\left(at   + s \sqrt{1-a^2}\sqrt{1-t^2}\right){\rm d}\mu^{(\gamma-1/2)}\ .
\end{equation}
 We refer to the  $K_a$ as the {\em correlation operators}
 for reasons that will be explained in Section~\ref{Bl}.
\end{defi}

With this definition and (\ref{gegen2}), (\ref{boc12}) can be written as
\begin{equation}\label{boc11b}
K_a p^{(\gamma)}_n(x)  = \frac{  p^{(\gamma)}_n(a)}{p^{(\gamma)}_n(1)}p^{(\gamma)}_n(x)\ .
\end{equation}
Thus $\{p^{(\gamma)}_n(a)/p^{(\gamma)}_n(1) \}_{n\ge 0}$ is the eigenvalue sequence of $K_a$,
and since the eigenvalues are real and the eigenfunctions are orthogonal, it follows that $K_a$ is self adjoint, as noted in general in the proof of Theorem~\ref{BHT}.

\noindent{\bf Proof of Gegenbauer's identity:} 
The starting point is a {\em direct} proof that $K_a$ is self adjoint.
 From (\ref{Kdef}) we find,
$$\langle K_a f, g\rangle_{L^2(\mu^{(\gamma)}) } = 
c_\gamma c_{\gamma-1/2}\int_{-1}^1 \int_{-1}^1g(t) 
f\left(at   + s \sqrt{1-a^2}\sqrt{1-t^2}\right)(1-s^2)^{\gamma-1}{\rm d}s(1-t^2)^{\gamma-1/2} {\rm d}t\ .$$
With the change of variables 
$u = at   + s \sqrt{1-a^2}\sqrt{1-t^2}$, the integral over $s$ becomes
$$= \int_{at-\sqrt{1-a^2}\sqrt{1-t^2}}^{at+\sqrt{1-a^2}\sqrt{1-t^2}}f(u)\frac{((1-a^2)  - (u^2+t^2 - 2atu))^{\gamma-1}}{(1-a^2)^{\gamma-1/2}}du,
$$
so that
\begin{equation}\label{genquad}
\langle K_a f, g\rangle_{L^2(\mu^{(\gamma)}) } = 
c_\gamma c_{\gamma-1/2} \int_{-1}^1 \int_{-1}^1g(t) 
f(u)\frac{((1-a^2)  - (u^2+t^2 - 2atu))_+^{\gamma-1}}{(1-a^2)^{\gamma-1/2}}{\rm d}u{\rm d}t\ ,
\end{equation}
where $(\cdot )_+$ denotes the positive part.  Thus, $K_a$ is self adjoint 
on $L^2(\mu^{(\gamma)})$. 

For the rest, we reurn to the original defining formula (\ref{Kdef}) First, by  the symmetry properties of
$\mu^{(\gamma-1/2)}$,  $K_a$ maps polynomials of degree $n$ to polynomials of
degree $n$. It follows that the spectrum is discrete and the eigenfunctions 
are polynomials that are orthogonal
with respect to the measure $\mu^{(\gamma)}$, and hence are the  
$p_n^{(\gamma)}$. Let $\lambda_n$ be the eigenvalue corresponding to 
$p_n^{(\gamma)}$; i.e.,
$\lambda_n p_n^{(\gamma)}(t) = K_a  p_n^{(\gamma)}(t)$.
Taking the limit $t\to 1$ on both sides, using the Dominated Convergence 
Theorem  yields,
\begin{equation}
\lambda_n p_n^{(\gamma)}(1) = p_n^{(\gamma)} (a)
\end{equation}
which  immediately gives Gegenbauer's identity in the form (\ref{boc11b}). 
\lanbox

In what follows, we shall make repeated use of the mechanism illustrated in 
our proof of Gegenbauer's identity, and the next theorem paves the way for
its broader application:

\begin{thm}[{\bf Evaluation Formula}] \label{dec2}  Let $X$ be a closed 
interval in $\R$, and let $\mu$ be a regular Borel probability measure
whose support is $X$. Let $\{f_n\}_{n\ge 0}$
be the unit orthonormal basis for 
$L^2(\mu)$ consisting of the normalized orthogonal polynomial for $\mu$.  
Suppose that for each $z\in X$,
$K_z$ is an operator on  $L^2(\mu)$ with the following properties:
\smallskip

\noindent{\it (1)} $K_z$ is self adjoint on $L^2(\mu)$.

\smallskip

\noindent{\it (2)} If $f$ is a polynomial of degree no greater than $n$, then so is $K_z f$. 
\smallskip

\noindent{\it (3)} There exists an $x_0\in X$ such that
for any continuous function $f$, and any $z\in X$,
\begin{equation}\label{lim}
\lim_{x\to x_0}K_z f(x) = f(z)\ .
\end{equation}

Then for each $n$, $f_n(x_0) \ne 0$, and  for each $x$ ${\displaystyle K_z f_n(x) = \frac{f_n(z)}{f_n(x_0)}f_n(x) }$, so that
if $K_z$ is a Markov operator, then 
$\{f_n(z)/f_n(x_0)\}_{n\geq 0}$ is a Markov sequence for  $\{f_n\}_{n\geq 0}$.
\end{thm}

\noindent{\bf Proof:} Properties {\it (1)} and  {\it (2)}  immediately imply that each $K_z$ is diagonalized by polynomials that are orthogonal in $L^2(\mu)$, so that the eigenfunctions of $K_z$ are the
$f_n$.  To determine the eigenvalues,  start from the definition of the $n$th eigenvalue $\lambda_n$, 
$K_z f_n(x) = \lambda_n f_n(x)$,
and take the limit $x\to x_0$. By {\it (3)} we obtain
$ f_n(z) = \lambda_n f_n(x_0)$,
which tells us $f_n(x_0) \ne 0$ and $\lambda_n =f_n(z)/f_n(x_0)$ . \lanbox

To summarize, through an analysis of the operators $K_a$, based on the three properties highlighted in the previous theorem, we obtain a self-contained proof of Gegenbauer's identity, and hence Bochner's Theorem. Is there a analogous family of operators that gives Gasper's theorem? 

\subsection{Gasper's Theorem}

For $h\in {\mathcal  C}([-1,1])$ define $K_{a,0}$ by,
\begin{equation} \label{s44}
(K_{a,0}h)(t) = 
\int_0^1 \int_0^\pi h\left[ a^2(1+t)-1 + b^2 (1-t) r^2 +2 abr(1-t^2)^{1/2} \cos \theta \right] dm_{\alpha, \beta}(r,\theta)
\end{equation}
where
\begin{equation} \label{dmalpha}
{\rm d}m_{\alpha,\beta}(r,\theta) =
\frac{2 \Gamma(\alpha+1)}{\sqrt \pi \Gamma(\beta+1/2) \Gamma(\alpha-\beta)} (1-r^2)^{\alpha-\beta-1}
r^{2\beta+1}\sin^{2\beta}\theta  {\rm d}r   {\rm d}\theta \ 
\end{equation}
is a probability measure.
We now have,

\begin{lm}\label{gaspa} For all $a\in (-1,1)$ and $\alpha>\beta>-1/2$, the operator
$K_{a,0}$  on ${\mathcal  C}([-1,1])$ as defined in (\ref{s44}) has the following properties:

\noindent{\it (1)} $K_{a,0}$ is self adjoint on $L^2(\mu^{\alpha,\beta})$.

\smallskip

\noindent{\it (2)} The space of polynomials of any fixed degree is invariant under $K_{a,0}$.

\smallskip

\noindent{\it (3)} For any continuous function $h$,
${\displaystyle \lim_{t\to 1} K_{a,0} h(t ) = h(2a^2-1)}$.
\end{lm}

\smallskip

\noindent{\bf Proof:} Given the explicit formula (\ref{s44}), the proof of {\it (2)} follows from the form of $\mu^{\alpha,\beta}$ which shows that only even powers of $\cos\theta$ are nonzero when integrating over $\theta$. Part {\it {3}} follows from  the evaluation property (\ref{lim}) and the Dominated Convergence Theorem. It is only {\it (1)} that requires more work.

We now use a sequence of variable changes due to Koornwinder \cite{korn}, but
for a different purpose. We shall contrast our use of it with Koornwinder's in the final section of the paper, but for now, suffice it to say that Koornwinder was not concerned with self-adjointness, which is the issue before us.

Consider $h_1$ and $h_2$ in ${\mathcal  C}([-1,1])$  Then  by  (\ref{s44}) and the change of variables $t =2s^2-1$, 
$\langle h_1, K_{a,0} h_2\rangle_{L^2(\mu_{\alpha,\beta})}$ is a constant multiple
of
\begin{align*}
q(h_1,h_2)   &:= \int^1_0\int^1_0\int^\pi_0 h_1(2s^2-1)h_2((2a^2s^2-1)+ 2b^2(1-s^2)r^2+
4abrs\sqrt{1-s^2}\, \cos\theta)\\
&\times \quad\sin^{2\beta} \theta (1-r^2)^{\alpha-\beta-1}r^{2\beta+1}(1-s^2)^{\alpha}
s^{2\beta+1}{\rm d}s\, {\rm d}r\, {\rm d}\theta .\end{align*}
We must show that $q(h_1,h_2) = q(h_2,h_1)$. 

The first step  is to  replace $(1-r^2)^{\alpha - \beta -1}$ by 
$(1-r^2)_+^{\alpha - \beta -1}$, and extend the domain of integration in $r$ to  $(0,\infty)$. 
The point is that we may then
regard the integration over $r$ and $\theta$ as an integration over the upper half plane
in $\R^2$.  Changing to Cartesian coordinates $x$ and $y$
yields
\begin{align*}
q(h_1,h_2)&=\int^1_0 h_1(2s^2-1)(1-s^2)^{\alpha}s^{2\beta+1}\\
& \times\left(\int^\infty_{-\infty}\int^\infty_0
h_2(2(cx+as)^2+2c^2y^2-1)(1-x^2-y^2)^{\alpha-\beta-1}_+
y^{2\beta}dy\, dx\right){\rm d}s\end{align*}
where $c=b(1-s^2)^{1/2}$.

The second step is to translate and scale, making the change of variables $y'=cy$ and $x' =cx+as$. 
Then since
$(1-s^2)^\alpha y^{2\beta}{\rm d}x {\rm d}y = b^{-2\alpha} c^{2(\alpha - \beta -1)}{\rm d}x' {\rm d}y'$,
this  yields
\begin{align*}
q(h_1,h_2)  &=\int^1_0 h_1(2s^2-1)s^{2\beta+1} b^{-2\alpha}\\
&\times \left(\int^\infty_{-\infty} \int^\infty_0
h_2(2({x'}^2+{y'}^2)-1) (c^2-(x' -as)^2-{
  y'}^2)^{\alpha-\beta-1}_+{y'}^{2\beta}{\rm d} y'\, {\rm d}x' 
\right){\rm d}s\ .
\end{align*}

Finally, the third step is to  change back to polar coordinates; i.e., make
the change of variables $(x', y')\to (\rho,\phi)$. This yields, making crucial use of $a^2+b^2=1$, 
\begin{align*}
q(h_1,h_2)  &=\int^1_0 h_1(2s^2-1)s^{2\beta+1} b^{-2\alpha}\\
&\times\left(\int^\infty_{0} \int^{\pi}_0
h_2(2\rho^2-1) (b^2-s^2-\rho^2 +2a\rho s\cos{\phi})^{\alpha-\beta-1}_+ r^{2\beta+1}\,
{\rm d}  \phi\, {\rm d}  \rho \,\right) {\rm d} s\\
&=b^{-2\alpha}\int^1_0\int^1_{0} \int^{\pi}_0
h_1(2s^2-1)h_2(2\rho^2-1)\\
&\times (b^2-s^2-\rho^2 +2a\rho s\cos{\phi})_+^{\alpha-\beta-1} \rho^{2\beta+1}\,
s^{2\beta+1}\sin^{2\beta}\phi{\rm d} \phi\, {\rm d}  \rho\, {\rm d} s\, ,
\end{align*}
which finally renders  the symmetry manifest. \lanbox

\smallskip
\noindent{\bf Proof of Theorem~\ref{gasp1}} The case $\alpha=\beta>-1/2$ are contained in Bochner's theorem. 
For  $\alpha>\beta >-1/2$, Lemma~\ref{gaspa}  
implies that the family of operators $K_{a,0}$, $a\in [-1,1]$ satisfies the conditions of Theorem~\ref{dec2}
for the unit orthonormal basis $\{p_n^{(\alpha,\beta)}\}_{n\geq 0}$.    Then from the conclusion
of Theorem~\ref{dec2}, we may apply Theorem~\ref{BHT} to obtain Gasper's Theorem in this case.
The case $\alpha>\beta, \beta=-1/2$ follows since $p^{(\alpha,\alpha)}_{2n}(x)=p^{(\alpha,-1/2)}_n(2x^2-1)$ 
\lanbox

The remaining mystery at this point is where the operators $K_a$ and $K_{a,0}$ came from. 
In fact, the operator $K_a$ arose naturally in the work \cite{CCL2} on the Kac model  \cite{K}. The Kac model is a model from mathematical physics
for the trend to equilibrium in a gas of $N$ molecules interacting through binary collisions.  An analysis made in \cite{CCL2} of how the rate of equlibriation depends on $N$ for one dimensional velocities
reduces this issue to the determination of the eigenvalues of the operators $K_a$ (\ref{Kdef}), which measure correlations between the different particle's velocities. In the analysis of the Kac model for three dimensional velocities the following
operator naturally arises,

\begin{equation}\label{kackern}
{\mathcal  K}_a f(v) = \int_B f\left(av + \sqrt{1-a^2}\sqrt{1-|v|^2}y\right){\rm d}\nu_{m,N-1}(y)\  ,
\end{equation}
where
\begin{equation}\label{nuform}
{\rm d}\nu_{m,N}(v) = \frac{|S^{m(N-1) -1}|}{|S^{mN-1}|}(1 - |v|^2)^{(m(N-1)-2)/2}{\rm d}v\  ,
\end{equation}
and $B$ is the unit ball in $\R^m$, $|S^{d-1}|$ is the surface area of the unit sphere in $\R^d$. The number $a$ ranges from $-1$ to $1$, $N>2$ and $m> 1$.

Note the similarity of ${\mathcal  K}_a$ to the operator $K_a$ defined in (\ref{Kdef}). Of course ${\mathcal  K}_a$
acts on functions on the unit ball $B$, however  there is a direct connection to operators that act on functions
of $[-1,1]$, such as $K_{a,0}$. This follows from the fact that ${\mathcal  K}_a$ commutes with rotations and therefore preserves the class of
radial functions.  For $h\in {\mathcal  C}([-1,1])$, define $K_{a,0} h$ by
\begin{equation}\label{Kzero}
(K_{a,0} h)(2|v|^2-1) := ({\mathcal  K}_af)(v)\qquad {\rm where}\qquad f(v) := h(2|v|^2 -1)\ .
\end{equation}
This operator is well defined since ${\mathcal  K}_a$ preserves the class of radial functions.
A calculation, which we shall make in Section~\ref{GS5},  shows that for
\begin{equation}\label{ab}
\alpha =  (m(N-2)-2)/2 \ \ \ \ {\rm and} \ \ \  \ \beta = (m-2)/2 \
,\end{equation} 
the operators $K_{a,0}$ defined in (\ref{s44}), and (\ref{Kzero}) are the same. 
Thus, at least for the half integral values of $\alpha$ and $\beta$ in (\ref{ab}),
the apparently more complicated operator defined in (\ref{s44}) does indeed come from an operator bearing a striking resemblance to the one in Gegenbauer's identity.  Moreover, while $m$ and $N$ are integers in (\ref{kackern}), once
the radial part has been rewritten in the form (\ref{s44}), there is no obstacle to letting $\alpha$ and $\beta$
vary continuously.

The remarkable thing about this construction of $K_{a,0}$ is that it only uses one invariant subspace of the
operators ${\mathcal  K}_a$ to 
recover the known results of Gasper and Koornwinder.  However, there are other
invariant suspaces with a direct connection to Jacobi polynomials.
In fact, we shall see that for each integer $\ell>0$ there is a family of operators $K_{a,\ell}$
to which we may apply Theorem~\ref{dec2}.. 
For example, the following result
about Jacobi polynomials  is the analog of Gegenbauer's product formula (\ref{boc12}):

\begin{thm}\label{betterprod}For all $\alpha> \beta > -1/2$, and all non negative integers $\ell$,
\begin{eqnarray}\label{atl5}
 & &a^\ell \frac{p_n^{(\alpha,\beta+\ell)}(t)}{p_n^{(\alpha,\beta+\ell)}(1)}
p_n^{(\alpha,\beta+\ell)}(2a^2-1)\nonumber\\
 &=&
\int_0^1\int_0^\pi  p^{\alpha,\beta+\ell}_n\left(\left[ a^2(1+t) + b^2(1-t) r^2 + 2ab\sqrt{1-t^2} r\cos\theta    \right]-1\right)\nonumber\\
&\times &\left[\sum_{j=0}^\ell \left( \begin{array}{c} \ell \\
j 
\end{array}\right) a^{\ell-j}(br)^j \left({\frac{1-t}{1+t}}\right)^{j/2}  P_j^{(\beta)}(\cos \theta)\right]
{\rm d}m_{\alpha,\beta}(r,\theta)\ ,\nonumber\\
\end{eqnarray}
where $b = \sqrt{1-a^2}$ as before.
\end{thm}
Note that Gasper's formula appears for the case $\ell=0$.
As a consequence of this theorem we have
\begin{thm}\label{genko}  Consider any $\alpha> \beta> -1/2$ and any integer $\ell\ge0$. Let 
$P_\ell^{(\beta)}$ is the ultraspherical polynomial with the  normalization $P_\ell^{(\beta)}(1) =1$. Then for all $t\in [-1,1]$,
\begin{align*}
 \frac{p_n^{(\alpha,\beta+\ell)}(t)}{p_n^{(\alpha,\beta+\ell)}(1)}
=
\int_0^1\int_0^\pi  &\left[\frac{(1+t) - (1-t)r^2}{2} + i\sqrt{1-t^2}r\cos\theta \right]^{n}\\
\times &\left[\sum_{k = 0}^{\ell} 
 \left( \begin{array}{c} \ell \\
k 
\end{array}\right)
\left(\frac{1-t}{1+t}\right)^{k/2} (ir)^k P_k^{(\beta)}(\cos \theta)\right]
{\rm d}m_{\alpha,\beta}(r,\theta)\ .
\end{align*}
where
 ${\rm d}m_{\alpha,\beta}$ is given by (\ref{dmalpha}).
\end{thm}

The case $\ell = 0$ is a well--known formula of Koornwinder
\cite{Kor}. It was pointed out to us by an anonymous referee of a previous version of this work, that for $\ell>0$, the formula of Theorem \ref{betterprod} is equivalent (though not in such an obvious way)
to a product formula due to  Koornwinder and Schwartz 
for orthogonal polynomials on the so-called {\em parabolic biangle}; their equivalent formula is
(3.13) in \cite{KooSchwa}.

We state and indicate a direct proof, using our methods, of the  product formula for the parabolic biangle in the Appendix \ref{parabiangle}. In this Appendix we also discuss
the product formula on the triangle.

The rest of the paper is organized as follows:
In Section~\ref{Bl} we give a geometric and probabilistic interpretation of the operators $K_a$ and ${\mathcal  K}_a$. This shall
explain  our reasons for referring to them as  ``correlation operators''.  
It also yields a simple proof of their self adjointness, at least for the ``geometric'' values of $\gamma$,
 $\alpha$ and $\beta$. We then use  the rotational invariance of the operators ${\mathcal  K}_a$ to determine a sequence of
invariant subspaces for them, indexed by the non-negative integer $\ell$, and we study the spectrum of the restrictions
$K_{a,\ell}$ of ${\mathcal  K}_a$ to these invariant subspaces. Though for $\ell>0$, $K_{a,\ell}$ is not Markov, Theorem~\ref{dec2} is applicable nonetheless: The eigenvalues are again expressible in terms of ratios of Jacobi polynomials, and in this enable us to easily prove Theorems~\ref{betterprod} and \ref{genko}  for the ``geometric'' values of $\gamma$,
 $\alpha$ and $\beta$. 
 
 Then,  in Section~\ref{sect3},  we show how  $\gamma$,
 $\alpha$ and $\beta$ may be allowed to vary continuously, and thus prove 
 Theorems~\ref{betterprod} and \ref{genko} in full generality.

In Section~\ref{GSP6} we use the Laplace formula for ultraspherical polynomials and Theorem~(\ref{genko}) to obtain sharp bounds on ratios of Jacobi polynomials.  That is, we obtain sharp bounds on the eigenvalues of the extremal Markov operators, and these bounds give sharp information on the operator trace classes to which the extremal Markov operators belong. This  information is  then used to discuss the pointwise convergence  properties
of the  eigenfunction expansions for the  kernels assoicated with the operators in Bochner's and Gasper's theorem. 

In Section~\ref{h}, we discuss the history of Bochner's and Gasper's results and finally in the appendix we state and outline a proof of the
parabolic biangle and triangle polynomial product formula of \cite{KooSchwa} along the lines outlined in this section.

{\bf Acknowledgment:} The authors would like to thank Dominique Bakry
for illuminating discussions about his papers with Huet and Mazet, and for suggesting that ideas arising in
our work on the spectral gap for the Kac model with three dimensional velocities 
might lead to a self--contained proof of the Markov 
sequence problem for Jacobi polynomials. We also thank an anonymous referee of a previous version of our
paper for pointing out the connection of Theorem \ref{betterprod} and the parabolic biangle polynomial product formula
of Koornwinder and Schwartz \cite{KooSchwa}.



\section{The geometric cases}\label{Bl}

The proof of the product formula for Jacobi polynomials with
$\alpha =  (m(N-2)-2)/2$ and $ \beta = (m-2)/2$, where $m$ and $N$ 
are  positive integers as in (\ref{ab}), is particularly simple because of a geometric picture
for the correlation operator ${\mathcal  K}_a$ in these cases. In this section, we shall present a complete
proof for these cases, which we call the geometric cases. Then in the next section, we  shall complete
our analysis by showing that while the geometric picture only makes sense for integer values of $m$ and $N$, certain  formulas and results that one derives using the geometric picture retain their validity as $m$ and $N$ are allowed to vary continuously.

As we have noted, the operators $K_a$ and ${\mathcal  K}_a$ arose in the study of the Kac model, where they measured correlations. We start by explaining the geoemtry behind the simple operator $K_a$, for which the geometric vales of  $\gamma = (N-2)/2$, $N$ a positive integer. 

\subsection{The geometric origins of the correlation operator $K_a$}

As is well known, when $\gamma = (N-2)/2$, $\mu^{(\gamma)}$ is simply the image of the uniform probability measure $\sigma^N$ on $S^{N-1}$, the unit sphere in $\R^N$, under the map $x\mapsto x\cdot \widehat e$, where $\widehat e$ is any unit vector in $\R^N$. That is, if $\widehat e$ is any unit vector in $\R^N$, and $f$ is any bounded measurable function on $[-1,1]$, then
$$\int_{S^{N-1}} f(x\cdot \widehat  e){\rm d}\sigma_{N} = \int_{-1}^1f(t){\rm d}\mu^{((N-2)/2)}(t)\ .$$

Let $\widehat u_1$ and $\widehat u_2$ be any  
two unit vectors in $\R^N$, and define 
 a bilinear form  $q_{\widehat u_1,\widehat u_2}$
 on $L^2(\mu^{((N-2)/2)})$ by
\begin{equation}\label{boc21}
q_{\widehat u_1,\widehat u_2}(f,g) = 
\int_{S^{N-1}}f(x\cdot \widehat u_1)g(x\cdot \widehat u_2){\rm d}\sigma_{N}\ .
\end{equation} 

We claim
that $q_{\widehat u_1,\widehat u_2}(f,g)$ is symmetric in $f$ and $g$, and
depends on the choice of 
$\widehat u_1$ and $\widehat u_2$ only through $a:= \widehat u_1\cdot \widehat u_2$ 

To see this, 
let $T$ be the reflection in $\R^N$ about the hyperplane orthogonal to $\widehat u_2 - \widehat u_1$. Then
$T(\widehat u_2) = \widehat u_1$ and  $T(\widehat u_1) = \widehat u_2$, and hence, by the invariance of ${\rm d}\sigma_{N}$ under orthogonal transformation of $\R^N$,
$$q_{\widehat u_1,\widehat u_2}(f,g) = q_{\widehat u_2,\widehat u_1}(f,g) = q_{\widehat u_1,\widehat u_2}(g,f)\ .$$

A similar argument using a rotation that fixes, say,  $\widehat u_2$ shows that 
this bilinear form depends on $\widehat u_2$ and $\widehat u_1$ only 
through $a:=\widehat u_1\cdot \widehat u_2$, and this established the claim.

We may now use the quadratic form $q_{\widehat u_1,\widehat u_2}(f,g)$ to define an operator $K_a$ where
$a =  \widehat u_1\cdot \widehat u_2$. 
It turns out that
this operator associated is exactly the operator $K_a$ defined in (\ref{Kdef}):

\begin{prop}\label{sp1} For any $N>1$ and any $-1 < a < 1$,  and all continuous functions $f$ and $g$ on $[-1,1]$,
\begin{equation}\label{quadform}
\langle K_a f, g\rangle_{L^2(\mu^{((N-2)/2)})} = 
\int_{S^{N-1}} f(x\cdot \widehat u_1)g(x\cdot \widehat u_2){\rm d}\sigma_{N} \ .
\end{equation}
where $K_a$ is the operator on $L^2(\mu^{(\gamma)})$ for $\gamma = (N-2)/2$
defined in (\ref{Kdef}).
\end{prop}

\noindent{\bf Proof:} This is a calculation based on he following system of coordinates on $S^{N-1}$:
Define 
$$\phi: S^{N-2}\times [-1,1] \to S^{N-1}$$
by 
$\phi(y,t) = (\sqrt{1-t^2}y_1\dots,\sqrt{1-t^2}y_{N-1}, t)$.
Evidently for any $y\in S^{N-2}$ and any $t\in [-1,1]$, $\phi(y,t) \in S^{N-1}$.
It is then easy to check, as in \cite{CCL2},  that for any function $h$ on $S^{N-1}$,
$$\int_{S^{N-1}} h(x){\rm d}\sigma_{N}(x) = 
\int_{-1}^1\left[\int_{S^{N-2}} h(\phi(y,t)){\rm d}\sigma_{N-1}(y)\right]{\rm d}\mu^{((N-3)/2)}(t)\ .$$
We now apply this to the integral in (\ref{quadform}). Let $\{\widehat e_1,\dots,\widehat e_N\}$ be the standard
orthonormal basis for $\R^N$. Take
\begin{equation}\label{u1u2def}
\widehat u_1 = \widehat e_N \qquad{\rm and}\qquad \widehat u_2 = a \widehat e_N + \sqrt{1-a^2}\widehat e_{N-1}\ .\end{equation}
Then with  $h(x) = f(x\cdot \widehat u_1)g(x\cdot \widehat u_2)$, we obtain (\ref{quadform}). \lanbox

\medskip

It is clear from (\ref{quadform}) that $K_a$ is self adjoint for $\gamma = (N-2)/2$.  Once one knows the self-adjointness  for these special values of $\gamma$,
it is natural to seek a direct proof -- without lifting the functions onto spheres. What one finds is the
``intrinsic'' quadratic form representation (\ref{genquad}) that we gave 
in the introduction when we proved Gegenbauer's formula. While (\ref{genquad}) may be less elegant than (\ref{quadform}), it has the advantage that it is valid for all $\gamma > 1/2$.

We close this subsection by giving
simple  probabilistic interpretation $K_a$ which explain our use of the term ``correlation'':
Think of  $S^{N-1}$, equipped with ${\rm d}\sigma_N$ as a probability space, and think of
$f(x\cdot \widehat u_1)$ as a random variable on this probability space.  Then, the conditional
expectation of $f(x\cdot \widehat u_1)$ given $x\cdot \widehat u_2$ is the function $h(x\cdot \widehat u_2)$
such that 
${\rm E}\left[h(x\cdot \widehat u_2)g(x\cdot \widehat u_2)\right] = {\rm E}\left[f(x\cdot \widehat u_1)g(x\cdot \widehat u_2)\right]$
for all continuous bounded functions $g$. Since 
$$\langle K_a f, g\rangle_{L^2(\mu^{((N-2)/2)})}  = 
{\rm E}\left[K_af(x\cdot \widehat e)g(x\cdot \widehat e)\right] $$
for any unit vector $\widehat e$, in particular for $\widehat e = \widehat u_2$, we see from
(\ref{quadform})
that
\begin{equation}\label{condi}
K_a f(t) = {\rm E}\{ f(x\cdot \widehat  e_2) \ |\ x\cdot \widehat  e_1 = t\ \}\ .
\end{equation}

\if false

\begin{remark}
Note that the bilinear form $\langle f, {\mathcal K}_a g\rangle$ is exactly the convolution of the measures
$f(x)dx, g(x)dx$ on $[-1,1]$ in the sense of (\ref{convol}), for the special case where $\gamma = (N-2)/2$.
\end{remark}

We now explain how this conditioning operation can produce a family of operators satisfying the key property 
 {\it (3)} of Theorem \ref{dec2}.

\begin{prop}\label{sp2} Let $(\Omega,{\mathcal  S},\mu)$ be a probability space on which there are random variables
$X_a$, $a\in [-1,1]$ such that
$${\rm ess\ sup}(X_a) =1\qquad{\rm for \ all} \quad a \in [-1,1]$$
and with the property that for every $\epsilon>0$ there is a $\delta>0$
such that
$$X_1\ge 1-\delta \Rightarrow |X_a -a| \le \epsilon\ .$$
Then for any continuous and bounded function $f$ on $[-1,1]$,
$$\lim_{t\to 1} {\rm E}\{ f(X_a) \ |\ X_1 =  t\ \} = f(a)\ .$$
\end{prop}

\noindent{\bf Proof:} Fix $a$ and $f$.  Then by the continuity of $f$,  for any $\gamma>0$, 
there is  $\epsilon>0$ so that $|b-a| \le \epsilon \Rightarrow |f(b) - f(a)|\le \gamma$. By hypyoehtesis, 
there is a $\delta>0$ so that on the non-empty set $\{X_1 > 1-\delta\}$,  $|X_a -a| \le \epsilon$. Hence,
on the set $\{X_1 > 1-\delta\}$, $|f(X_a) - f(a)| \le \gamma$. Since $\gamma>0$ is arbitrary, the conclusion now follows.
\lanbox

This lemma may be applied in the case at hand with $(\Omega, {\mathcal  S}, \mu) = (S^{N-1},{\mathcal  B},\sigma_N)$
and $X_a = x\cdot (a \widehat e_N + \sqrt{1-a^2}\widehat e_{N-1})$ since then it is readily checked that
$$X_1 > 1 - \epsilon^2/2 \Rightarrow |X_a-a| < \epsilon\ .$$

This may shed some light on the origins of the evaluation property  {\it (3)} of Theorem \ref{dec2}, and of course the 
self-adjointness property {\it (1)} follows directly from the construction. however, we do not claim to have any sort of general prescription for writing down operators satisfying the conditions of  Theorem \ref{dec2}.

\fi
 
 \if false
As we have discussed in the introduction, associated to any product formula (\ref{product}) is a convolution operation. The product formula for the ultraspherical polynomials $\{P_n^{(\gamma)}\}_{n\ge 0}$; i.e., Gegenbauer's identity, has an interesting geometric interpretation when $\gamma = (N-2)/2$ 
for an integer $N\ge 2$.  

Fix some unit vector $e\in \R^N$. A measure on $S^{N-1}$ is {\em zonal about} $e$ if it is invariant under rotations preserving $e$.  Given any (finite, positive) measure $\lambda$ on $[-1,1]$, there is a unique
measure $\widehat \lambda$ on $S^{N-1}$ such that for all $\varphi\in {\mathcal C}([-1,1])$,
$$\int_{S^{N-1}}\varphi(u\cdot e){\rm d}\sigma_N(u) = \int_{[-1,1]}\varphi(t){\rm d}\lambda(t)\ .$$
Note that if $\lambda$ is absolutely continuous with respect to $\mu^{((N-2)/2))}$ with Radon-Nikodym derivative $h(t)$, so that ${\rm d}\lambda(t)  = h(t){\rm d}\mu^{((N-2)/2))}(t)$.
Then 
$${\rm d}\widehat \lambda(u) = h(u\cdot e){\rm d}\sigma_N(u)\ .$$
Let $\nu$ be another (positive, finite) measure on $[-1,1]$ that is absolutely continuous with respect to $\mu^{((N-2)/2))}$ with Radon-Nikodym derivative $g(t)$ so that 
${\rm d}\widehat \nu(u) = g(u\cdot e){\rm d}\sigma_N(u)$. Define a positive linear functional on
${\mathcal C}([-1,1])$ by
${\displaystyle \varphi \mapsto \int_{S^{N-1}}\int_{S^{N-1}}\varphi(u\cdot v){\rm d}\widehat \lambda(u) 
{\rm d}\widehat \nu(v)}$.
By the Riesz-Markov Theorem, there is a measure $\lambda\star\nu$ on $[-1,1]$ so that
\begin{equation}\label{prfo3}
\int_{[-1,1]}\varphi(t){\rm d}(\lambda\star\nu) = 
\int_{S^{N-1}}\int_{S^{N-1}}\varphi(u\cdot v){\rm d}\widehat \lambda(u) 
{\rm d}\widehat \nu(v)\ .
\end{equation}

\begin{thm}
The convolution operation defined by (\ref{prfo3}) for $\gamma = (N-2)/2$, $N\ge 2$ an integer,
coincides with the one defined by Gegenbauer's product formula for ultraspherical polynomials in the abstract formula (\ref{product}).
\end{thm}

\noindent{\bf Proof:}
It suffices to consider  the absolutely continuous case with ${\rm d}\widehat \lambda(u) = h(u\cdot e){\rm d}\sigma_N(u)$
and ${\rm d}\widehat \nu(v) = g(v\cdot e){\rm d}\sigma_N(v)$, the right hand side of 
(\ref{prfo3}) becomes
$$
\int_{S^{N-1}}\int_{S^{N-1}}\varphi(u\cdot v) g(v\cdot e) h(u\cdot e){\rm d}\sigma_N(u){\rm d}\sigma_N(v)\ .
$$
Since this is independent of $e$, we may integrate over $e$, and by Fubini's Theorem, we
may do the integration over $e$ first, so  that the  right hand side of 
(\ref{prfo3}) becomes
$$
\int_{S^{N-1}}\int_{S^{N-1}}\varphi(u\cdot v)\left(
\int_{S^{N-1}}g(v\cdot e) h(u\cdot e){\rm d}\sigma_N(e)\right){\rm d}\sigma_N(u){\rm d}\sigma_N(v)\ .
$$
Now let $k_z(x,y)$ denote the kernel for the operator $K_z$, noting that $K_z$ does indeed posses such a kernel whose explicit form can be read off from (\ref{genquad}).  We have
$$ \langle K_{z} g,h\rangle_{L^2(\mu^{((N-2)/2)})} = \int_{[-1,1]}\int_{[-1,1]} k_z(x,y){\rm d}\lambda(x){\rm d}\nu(y)\ .$$
Then  by (\ref{quadform}), 
$$
\int_{S^{N-1}}g(v\cdot e) h(u\cdot e){\rm d}\sigma_N(e) =
\int_{[-1,1]}\int_{[-1,1]} k_z(x,y){\rm d}\lambda(x){\rm d}\nu(y)\ ,$$
and hence
\begin{eqnarray}\int_{[-1,1]}\varphi(t){\rm d}(\lambda\star\nu) &=&
\int_{S^{N-1}}\int_{S^{N-1}}\varphi(u\cdot v)\left(
\langle K_{u\cdot v} g,h\rangle_{L^2(\mu^{((N-2)/2)})} \right){\rm d}\sigma_N(u){\rm d}\sigma_N(v)
\nonumber\\
&=&
\int_{[-1,1]}\varphi(z) \left(
\int_{[-1,1]}\int_{[-1,1]} k_z(x,y){\rm d}\lambda(x){\rm d}\nu(y)
\right){\rm d}\mu^{((N-2)/2)}(z)\nonumber
\end{eqnarray}
Now comparing with (\ref{prfo1b}) and (\ref{convol}), we see that the convolution defined by (\ref{prfo3})
coincides with the one defined using Gegenbauer's product formula.  \lanbox
\medskip

\noindent{There} is a different geometric picture for this convolution explained in Section 5.3 of \cite{BH}.

\fi


\subsection{The geometric origins of the correlation operator ${\mathcal  K}_a$}

The Jacobi polynomial version of Proposition~\ref{sp1}, leading to ${\mathcal  K}_a$ instead of $K_a$,  is only slightly more complicated than the original.
Note that any vector $x\in \R^{mN}$
can be written as an $N$--tuple of vectors in $\R^m$,
$x = (x_1,\dots,x_N)$, and hence may be identified with the $m\times N$ matrix
\begin{equation}\label{boc51}
[x] = [x_1,\dots,x_N]
\end{equation}
whose $j$th column is $x_j$. Then for any vector $u\in \R^N$, the matrix product $[x]u$ is well defined in
$\R^m$. It is easy to see that if $x\in S^{mN -1}$ and $\widehat u\in S^{N-1}$,
 then $[x]\widehat u$ lies in $B$, the unit ball in $\R^m$. 
 
Therefore, given two unit vectors $\widehat u_1$ and $\widehat u_2$ in $\R^N$, and
any two functions $f$ and $g$ on $B$, define 
\begin{equation}\label{boc52}
q_{\widehat u_1,\widehat u_2}(f,g) = 
\int_{S^{mN-1}}f ([x]\widehat u_1) g([x]\widehat u_2){\rm d}\sigma_{mN}\ .
\end{equation}
As before, this will depend only on the choices of  $\widehat u_1$ and $\widehat u_2$ through
$a = \widehat u_1\cdot \widehat u_2$. Hence we may use this bilinear form to
define a family of self-adjoint Markov operators on $L^2({\rm d}\nu_{m,N})$. Our next proposition says that
the operators we obtain this way are exactly the ${\mathcal  K}_a$:

\begin{prop}\label{dec4}
For any $N>2$ and $m>1$, and any $-1< a<1$, and all $f,g\in {\mathcal  C}(B)$, 
\begin{equation}\label{quadformmn}
\langle f, {\mathcal  K}_a g\rangle_{L^2(B,\nu_{m,N})} = 
\int_{S^{mN-1}}f ([x]\widehat u_1) g([x]\widehat u_2){\rm d}\sigma_{mN}\ ,
\end{equation}
where ${\mathcal  K}_a$ is the operator defined in (\ref{kackern}).
\end{prop}

\noindent{\bf Proof:}  We proceed exactly as in the proof of (\ref{quadform}).
Define 
$$\phi: S^{m(N-1)-1}\times B \to S^{mN-1}$$
by 
$\phi(y,v) = (\sqrt{1-|v|^2}y_1\dots,\sqrt{1-|v|^2}y_{N-1}, v)$.
It is then easy to check, as in \cite{CCL2},  that for any function $h$ on $S^{mN-1}$,
$$\int_{S^{mN-1}} h(x){\rm d}\sigma_{mN}(x) = 
\int_{B}\left[\int_{S^{m(N-1)-1}} h(\phi(y,v)){\rm d}\sigma_{m(N-1)}(y)\right]{\rm d}\nu_{m,N}(v)\ .$$
We now apply this to the integral in (\ref{boc52})
with $\widehat u_1$ and $\widehat u_2$ given by (\ref{u1u2def}). 
With  $h(x) = f([x]\widehat u_1)g([x]\widehat u2)$, we obtain (\ref{kackern}). \lanbox

As before, each ${\mathcal  K}_a$ is a self adjoint Markov operator on $L^2({\rm d}\nu_{m,N})$,
and has an interpretation as a conditional expectation operator:
${\mathcal  K}_a$, acting on functions on $B$, such that for all $v\in B$,
\begin{equation}\label{condi2}
{\mathcal  K}_a g(v) = {\rm E}\{ g([x]\widehat  e_2) \ |\ [x]\widehat  e_1 = v\ \}\ .
\end{equation}

In the next subsection, we exploit the self-adjointness of ${\mathcal  K}_a$ to obtain the 
the product formula for Jacobi polynomials in the geometric cases.





\subsection{Spectral analysis of ${\mathcal  K}_a$ and a product formula  in the geometric cases} \label{GS5}

In this subsection we study the operator ${\mathcal  K}_a$ restricted to various invariant subspaces. As we have seen, the restriction of 
${\mathcal  K}_a$ to the subspace of rotationally invariant subspaces gives Gasper's kernel $K_{a,0}$. The study of ${\mathcal  K}_a$ on other invariant subspaces leads to the Theorems~\ref{betterprod} and \ref{genko}.


\begin{lm}\label{dec66} For all $a\in (-1,1)$, and all $m>1$, $N>2$, ${\mathcal  K}_a$ has the following properties:
\smallskip

\noindent{\it (1)} ${\mathcal  K}_a$ is self adjoint on $L^2(\nu_{m,N})$.

\smallskip

\noindent{\it (2)} If $f$ is a polynomial of degree $n$ on $B$, then so is ${\mathcal  K}_af$. 
\smallskip

\noindent{\it (3)} For any continuous function $f$, and any unit vector $\widehat e$,
${\displaystyle
\lim_{t\to 1}{\mathcal  K}_a f(t\widehat e ) = f(a\widehat e)}$.

\smallskip

\noindent{\it (4)} For any rotation $R$ on $\R^m$,
${\displaystyle {\mathcal  K}_a(f\circ R) = ({\mathcal  K}_af)\circ R}$.
\end{lm}

\noindent{\bf Proof:}  We argue very much as we did in the ultraspherical case, except of course for the proof of {\it (4)}, which is a new multidimensional feature.

Propisition~\ref{dec4}, which expresses ${\mathcal K}_a$ in terms of a quadratic form
immediately yields  {\it (1)}.  As for {\it (2)},
note that  ${\rm d}\nu_{m,N-1}(s)$ is even in $s$. Therefore, if $m$ is any integer, all of the terms
that are of odd degree in $s$ that one obtains upon expansion of 
$\left(at   + s \sqrt{1-a^2}\sqrt{1-t^2}\right)^m$
drop out of the integral. Hence, what remains is a polynomial in $t$ of degree $m$.  

Further, {\it (3)} follows by the dominated convergence formula; take the limit under the integral sign, and use
$$\lim_{t\to 1}f\left(at\widehat e   + s \sqrt{1-a^2}\sqrt{1-t^2}\right) = f(a\widehat e)\ .$$
This is independent of $s$, and since $\nu_{m,N-1}$ is a probability measure, {\it (3)} now follows. 
Finally, {\it (4)} follows from the rotational  invariance of $\nu_{m,N-1}$. \lanbox

Since ${\mathcal  K}_a$ commutes with rotations we can study its action
on the irreducible subspaces of the rotation group.   We begin by considering the action of ${\mathcal  K}_a$
on the radial functions on $B$, and shall deduce an elegant product formula for Jacobi polynomials directly from Lemma~\ref{dec66}.

Note that if $q$ is a polynomial of degree at most $n$
in one real variable, and the function $f$ on $B$ is defined by $f(v) = q(|v|^2)$, then by parts 
{\it (2)} and {\it (4)} of Lemma ~\ref{dec66},
${\mathcal K}_af(v)$
is again of this same form -- a polynomial of degree at most $n$ in $|v|^2$. Thus,
the  subspace of such functions is invariant under ${\mathcal K}_a$. 

Since by part {\it (1)}  of  Lemma ~\ref{dec66}, ${\mathcal K}_a$ is self-adjoint
on $L^2(\nu_{m,N})$, it may be diagonalized on each invariant subspace. It easily follows from here that  for each integer
$n\ge 0$, there is a polynomial $q_n$ such that with $f_n(v) = q_n(|v|^2)$, $f_n$ is an eigenvector of ${\mathcal K}_a$
with eigenvalue $\lambda_n(a)$, and that the $f_n$, appropriately normalized constitute a unit orthonormal basis for the subspace of radial functions in $L^2(\nu_{m,N})$.
By the explicit form of $\nu_{m,N} $ given in (\ref{nuform}),
$$\int_B f_v(v)^2{\rm d}\nu_{m,N}  = 
C_{m,N}\int_{[0,1]}q_n^2(r^2)(1-r^2)^{m(N-2)-2)/2}r^{m-1}{\rm d}r$$
where
$C_{m,N}$ is a normalization constant.  Making the change of variables $t = r^2$, one now recognizes the $q_n$ as being the Jacobi polynomials for $\alpha,\beta$ given by
(\ref{ab}), translated and scaled so the domain is $[0,1]$ instead of $[-1,1]$.

To determine the eigenvalues $\lambda_n(a)$, apply the evaluation formula, part {\it (3)}  of  Lemma~\ref{dec66}, to see that
$\lambda_n(a)q_n(1) = \lim_{|v|\to 1}\lambda_n(a)q_n(|v|^2) =
{\mathcal K}_af_n(v) = q_n(a^2)$. That is:
$$
\lambda_n(a) = \frac{q_n(a^2)}{q_n(1)}\ ,
$$
Since the restriction of ${\mathcal K}_a$ to the radial functions on $B$ is clearly positivity preserving, and clearly preserves the constants, it follows that the $\lambda_n(a) := q_n(a^2)/q_n(1)$
are a Markov sequence for the $\{q_n\}_{n\ge 0}$. 

Thus, condition (1) of Theorem~\ref{BHT} is satisfied, and as a consequence of 
Theorem~\ref{BHT}, we have therefore solved the Markov sequence problem for the Jacobi polynomials, and have proved a product formula for them, in the geometric cases. 
One can of course undo the scaling and translation, and write this all out explicitly for the usual Jacobi polynomial defined on $[-1,1]$. The resultis, of course, Gasper's product formula. We shall do this, but first notice that there is more to be obtained from the
analysis of ${\mathcal K}_a$: So far, we have only considered the restriction of ${\mathcal K}_a$
to the radial functions. The spectral analysis of ${\mathcal K}_a$ on other invariant subspaces prvides additional formulas identifying ratios of Jacobi polynomials as eigenvalues of self adjoint operators. We shall use these formula (and their extension to general values of $\alpha$ and $\beta$)
to prove Theorems~\ref{betterprod} and ~\ref{genko}.

\subsection{The spectral analysis of ${\mathcal K}_a$ on non-radial functions}

For each integer $\ell\ge 0$, let ${\mathcal  H}_\ell$ denote the space of harmonic polynomials on $\R^m$
that are homogeneous of degree $\ell$. 
Restricted to $B$, the functions in ${\mathcal  H}_\ell$ constitute a closed subspace in 
$L^2(\nu_{m,N})$, which we again denote by  ${\mathcal  H}_\ell$.

For each $\ell$,  ${\mathcal  H}_\ell$ is an eigenspace of $ {\mathcal  K}_a$. In fact, 
 for each $H\in  {\mathcal  H}_\ell$,
 \begin{equation}\label{eigell}
 {\mathcal  K}_a H(v) = a^\ell H(v)\ .
 \end{equation}
 That is, the restriction of  ${\mathcal  K}_a$ to $H\in  {\mathcal  H}_\ell$ is $a^\ell$ times the identity.  One way to see this is to use the mean value property of harmonic functions and the formula (\ref{nuform}).
 Since the measure ${\rm d}\nu_{m,N-1}$ is radially symmetric, we see that  ${\mathcal  K}_aH(v)$ = 
$H(av)$, which, by the homogeneity, is $a^\ell H(v)$.

There is another more algebraic argument that tells us somewhat more:

\begin{lm}\label{har} The spectrum of ${\mathcal  K}_a$ is discrete, and its eigenfunctions are of the
form $g(|v|^2)H(v)$, where $g(|v|^2)$ is a polynomial in $|v|^2$ and
$H \in  {\mathcal  H}_\ell$ for some $\ell$. Moreover, if $g(|v|^2)H(v)$ is an eigenfunction, then, so is
$g(|v|^2)\tilde H(v)$, for {\em any} non zero  $\tilde H\in {\mathcal  H}_\ell$.
\end{lm}

\noindent{\bf Proof:}
By Lemma \ref{dec66} the operator ${\mathcal  K}_a$ leaves the space of polynomials of degree $n$ invariant for any $n$. Hence, by the Weierstrass theorem the eigenfunctions consist of polynomials.
Further, since  ${\mathcal  K}_a$ commutes with rotations, any eigenfunction must be of the form
\begin{equation}\label{form}
F(v)=f(|v|) {\mathcal  Y}^\ell\left(\frac{v}{|v|}\right)
\end{equation}
where ${\mathcal  Y}^\ell$ is a spherical harmonic, i.e.,
${\displaystyle
{\mathcal  Y}^\ell\left(\frac{v}{|v|}\right) = |v|^{-\ell} H_\ell(v)
}$
where $H_\ell(v)$ is a homogeneous harmonic polynomial of degree 
$\ell$. We have to show that $f(|v|)/|v|^\ell$ is a polynomial
in $v$, i.e., a polynomial of the variable $|v|^2$.

Since $F(v)$ is a polynomial of degree $n$ we can write it as
${\displaystyle
F(v) =  \sum_{m=0}^n q_m(v)
}$
where $q_m(v)$ is homogeneous of degree $m$. In turn, each of these
polynomials can be expanded in terms of homogeneous harmonic polynomials, i.e.,
${\displaystyle
q_m(v) = H_m(v) + |v|^2 H_{m-2}(v) + |v|^4 H_{m-4}(v) + \cdots  \ .
}$
This shows that
\begin{equation}
F(v) = \sum_{k=0}^n g_k(|v|^2) H_k(v) 
\end{equation}
for some polynomials $g_k$. The result follows from (\ref{form}) 
and the orthogonality  properties of the spherical harmonics. The final statement follows from 
Schur's Lemma since ${\mathcal  K}_a$ commutes with rotations and rotations act irreducibly on ${\mathcal  H}_\ell$.
\lanbox

 Now, since polynomials on $[0,1]$ are 
uniformly dense in ${\mathcal  C}([0,1])$, it follows from the Lemma (and the fact that ${\mathcal  K}_a$ is Markov)  that 
 for each  function $g\in {\mathcal  C}([0,1])$, and each  $H\in  {\mathcal  H}_\ell$,   and all $a\in (-1,1)$,
 there is a   $\widetilde g_a\in {\mathcal  C}([0,1])$ so that
\begin{equation}\label{eigell2}
{\mathcal  K}_a f(v) = \widetilde g_a(|v|^2)H(v) \quad{\rm where}\quad f(v) = g(|v|^2)H(v)
\end{equation}
with $H\in  {\mathcal  H}_\ell$ being the same on both sides.  The transformation $g \mapsto
 \widetilde g$ is clearly linear, and as one sees from the proof of Lemma~\ref{har}, independent of the choice of $H$. We now use this transformation to  generalize the definition of the operator in (\ref{s44}).

To make efficient contact with the theory of Jacobi polynomials, it is better to write
our radial functions in the form $h(2|v|^2-1)$ instead of $g(|v|^2)$. For any non zero $H$ in any ${\mathcal  H}_\ell$,
we define ${\mathcal  V}_H$ to be the subspace of $L^2(\nu_{m,N-1})$ consisting of functions of the form
$$f(v) = h(2|v|^2 -1)H(v)\ ,$$
where $h$ is a function on $[-1,1]$. We then generalize the definition (\ref{s44}) as follows:

For each $\ell>0$, fix some non zero  $H \in {\mathcal  H}_\ell$.
  Then for $h\in {\mathcal  C}([-1,1])$, define $K_{a,\ell} h$ by
\begin{equation}\label{Kell}
(K_{a,\ell} h)(2|v|^2-1)H(v) = ({\mathcal  K}_af)(v)\qquad {\rm where}\qquad f(v) := h(2|v|^2 -1)H(v)\ .
\end{equation}
By the last statement in Lemma \ref{har}, $(K_{a,\ell} h)(2|v|^2-1)H(v)$ does not depend on the particular choice
of $H$ in ${\mathcal  H}_\ell$.
Further, by Lemma \ref{har} the  eigenfunctions of 
${\mathcal  K}_a$ are of the form
$$f_{n,\ell}(v) = h_{n,\ell}(2|v|^2 -1)H(v)$$
where $h_{n,\ell}$ is a polynomial of degree $n$. 
The eigenfunctions of $K_{a,\ell}$ are the polynomials $h_{n,\ell}$.  We now identify these as 
Jacobi polynomials:

First, we normalize our choice of $H\in {\mathcal  H}_\ell$ so that 
${\displaystyle
\int_{S^{m-1}}|H(v)|^2{\rm d}\sigma_m = 1
}$.
Then for any positive integers $n\ne p$,  integrating in polar coordinates and using (\ref{nuform})
 we find
\begin{align*}
0 &= \langle  f_{n,\ell},f_{p,\ell}\rangle_{L^2(\nu_{m,N})} =
\int_B h_{n,\ell}(2|v|^2 -1) h_{p,\ell}(2|v|^2 -1)|H(v)|^2{\rm d}\nu_{m,N}(v)\\
&= 
\frac{|S^{m(N-1) -1}|}{|S^{mN-1}|}\int_0^1h_{n,\ell}(2s^2-1)h_{p,\ell}(2s^2-1)
(1 - s^2)^{(m(N-1)-2)/2}s^{2\ell + m-1}{\rm d}s \ .\\
\end{align*}
Making the now familiar change of variables $t = 2s^2-1$, we find
$$\int_{-1}^1  h_{n,\ell}(t)h_{p,\ell}(t)
(1 - t)^{(m(N-2)-2)/2}(1+t)^{\ell + (m-2)/2}{\rm d}t = 0\ ,$$
which is the orthogonality relation defining the Jacobi polynomials $p_n^{(\alpha,\beta+\ell)}$ with $\alpha$ and $\beta$ given by (\ref{ab}).

We now determine the eigenvalues $\lambda_{n,\ell}(a)$ such that
$K_{a,\ell}h_{n,\ell} =  \lambda_{n,\ell}(a)h_{n,\ell}$.  By (\ref{Kell}), if we define
$f$ by $f(v) = h_{n,\ell}(2|v|^2-1)H(v)$, we have
$\lambda_{n,\ell}(a)f(v) = {\mathcal  K}_af(v)$.  Then, for any unit vector $\widehat e$ in $\R^m$ with
$H(\widehat e) \ne 0$, we have from part {\it (4)} of Lemma~\ref{dec66} that 
$\lambda_{n,\ell}(a)f(\widehat e) = f(a\widehat e)$.
which, by the homogeneity of $H$, means that $\lambda_{n,\ell}(a) = a^\ell h_{n,\ell}(2a^2-1)/h_{n,\ell}(1)$.
We summarize our conclusions in a lemma:

\begin{lm}\label{higherell} Fix dimensions $m>1$ and $N>2$, and let  $\alpha$ and $\beta$ be given by
(\ref{ab}).
Then for each integer $\ell\ge 0$ and each $a\in [-1,1]$, the operator $K_{a,\ell}$
is self adjoint on $L^2(\mu^{(\alpha,\beta+\ell)})$, and is diagonalized by the Jacobi polynomial basis
$\{p_n^{(\alpha,\beta+\ell)}\}_{n\ge 0}$. Moreover, the corresponding sequence of
eigenvalues $\{\lambda_{n,\ell}(a)\}_{n\ge 0}$ is given by 
\begin{equation}\label{lrat}
\lambda_{n,\ell}(a)= a^{\ell}\frac{p_n^{(\alpha,\beta+\ell)}(2a^2-1)}{p_n^{(\alpha,\beta+\ell)}(1)}\ .
\end{equation}
\end{lm}

\begin{remark}\label{indep}
 It is  worth remarking that the operators $K_{a,\ell}$ are {\em not} positivity preserving for $\ell> 0$. Nonetheless, the
eigenvalues  $\lambda_{n,\ell}(a)$ are eigenvalues of a Markov operator, namely ${\mathcal  K}_a$. 
\end{remark}

\section{General values of $\alpha$ and $\beta$}\label{sect3}

Our next goal is to extend this analysis of the previous section  to arbitrary values of $\alpha > \beta > -1/2$. We seek a direct expression of $K_{a,\ell}$, not explicitly involving ${\mathcal  K}_a$, so that we may then freely vary
the dimensions.  The following notation will be useful: For $v$ and $y$ in $B$, define
$$w(v,y,a) = av + \sqrt{1-a^2}\sqrt{1-|v|^2}y\ .$$
Then by Lemma \ref{dec4}, for any $f$ of the form $f(v) = h(2|v|^2-1)H(v)$, 
\begin{equation}\label{mon1a}
{\mathcal  K}_a f(v)  = \int_B h(2|w(v,y,a)|^2-1)H(w(v,y,a)){\rm d}\nu_{m,N-1}(y)\ .
\end{equation}

To proceed, we now make a judicious choice of $H$ to simplify the computations before us:

\noindent{$\bullet$} {\it We choose $H$ so that
$H(v/|v|)$ is the {zonal spherical harmonic} of degree $\ell$ with the axis along the unit vector
$\widehat e$ in $\R^m$. That is, }
\begin{equation}\label{mon2}
H(v) = |v|^\ell p_\ell^{((m-2)/2)}(\widehat e \cdot v/|v|)  = |v|^\ell p_\ell^{(\beta)}(\widehat e \cdot v/|v|)\ ,
\end{equation}
where, as before, $p_\ell^{(\beta)}$ is an ultraspherical polynomial, and $\beta$ is again given by (\ref{ab}).
 The advantage of our particular choice of $H$ is that $H(w(v,y,a))$ depends on 
$w(v,y,a)$ only through
$(w(v,y,a))\cdot \widehat e$ and $|w(v,y,a)|^2$. Specifically,
\begin{equation}\label{mon2a}
H(w(v,y,a)) =  |w(v,y,a)|^\ell p_\ell^{(\beta)}\left(\frac{(w(v,y,a))\cdot \widehat e}{|w(v,y,a)|}\right)\ .
\end{equation}
Note that
\begin{equation}\label{mon3}
w(v,y,a)\cdot \widehat e = as + \sqrt{1-a^2}\sqrt{1-s^2}r\cos\theta
\end{equation}
and
\begin{equation}\label{mon4}
|w(v,y,a)|^2   = a^2s^2 + (1-a^2)(1-s^2)r^2  + 2a\sqrt{1-a^2}\sqrt{1-s^2}rs\cos\theta
\end{equation}
where $s = |v|$, $r = |y|$, and $v\cdot y = sr\cos\theta$,
so that (\ref{mon2}) can be written as an integral over $r$ and $\theta$, 
using the measure defined
in (\ref{dmalpha}). These are the coordinates we used in the proof of Lemma~\ref{gaspa}  to ``liberate'' the values of $\alpha$ and $\beta$
in the $\ell = 0$ case.

By (\ref{Kell}) and (\ref{mon2}) evaluated at $v = s\widehat e$, we have
$$\left(K_{a,\ell}h\right)(2s^2-1)s^\ell H(\widehat e) = 
\int_B h(2|w(s\widehat e,y,a)|^2-1)H(w(s\widehat e,y,a)){\rm d}\nu_{m,N-1}(y)\ .$$
For the particular choice of $H$ made in (\ref{mon2}), this reduces to
\begin{equation}\label{atl1}
\left(K_{a,\ell}h\right)(2s^2-1) = s^{-\ell}
\int_B h(2|w(s\widehat e,y,a)|^2-1)|w(s\widehat e,y,a)|^{\ell}
P_\ell^{(\beta)}\left(\frac{(w(s\widehat e,y,a))\cdot \widehat e}{|w(s\widehat e,y,a)|}\right){\rm d}\nu_{m,N-1}(y)\ ,
\end{equation}
since $p_\ell^{(\beta)}(x)/ p_\ell^{(\beta)}(1) = P_\ell^{(\beta)}(x)$, the ultraspherical polynomial normalized by 
the condition $P_\ell^{(\beta)}(1) =1$.

Next, since the integrand depends on only on $s$, $r$ and $\cos\theta$, we can use (\ref{nuform}) and (\ref{dmalpha}) to write
this in terms of an integration against ${\rm d}m_{\alpha,\beta}$ with $\alpha$ and $\beta$ related to $m$ and $N$ through
(\ref{ab}):
\begin{equation}\label{atl2}
\left(K_{a,\ell}h\right)(2s^2-1) = s^{-\ell}
\int_0^1\int_0^\pi  h(2|w(s\widehat e,y,a)|^2-1)|w(s\widehat e,y,a)|^{\ell}
P_\ell^{(\beta)}\left(\frac{(w(s\widehat e,y,a))\cdot \widehat e}{|w(s\widehat e,y,a)|}\right){\rm d}m_{\alpha,\beta}(r,\theta)\ ,
\end{equation}

\begin{defi}\label{genkel} For all $\alpha> \beta > -1/2$, we define $K_{a,\ell}$ by formula (\ref{atl2}).
By the calculation just made, this coincides with the definition made in (\ref{Kell}) for 
$\alpha$ and $\beta$ satisfying (\ref{ab}).
\end{defi}

The next Lemma gives a more explicit formula for $K_{a,\ell}$.

\begin{lm}\label{better}For all $\alpha> \beta > -1/2$, and all non negative integers $\ell$,
\begin{eqnarray}\label{atl6}
K_{a,\ell}h(t) &=&
\int_0^1\int_0^\pi  h\left(\left[ a^2(1+t) + b^2(1-t) r^2 + 2ab\sqrt{1-t^2} r\cos\theta    \right]-1\right)\nonumber\\
&\times &\left[\sum_{j=0}^\ell \left( \begin{array}{c} \ell \\
j 
\end{array}\right) a^{\ell-j}(br)^j \left({\frac{1-t}{1+t}}\right)^{j/2}  P_j^{(\beta)}(\cos \theta)\right]
{\rm d}m_{\alpha,\beta}(r,\theta)\ ,\nonumber\\
\end{eqnarray}
where $b = \sqrt{1-a^2}$ as before.
\end{lm}

\noindent{\bf Proof:} First make the change of variables $t = 2s^2-1$ in the defining formula (\ref{atl2}). Under this change of variable,
(\ref{mon3})  and (\ref{mon4}) become
\begin{equation}\label{mon3b}
w(v,y,a)\cdot \widehat e = a\sqrt{\frac{1+t}{2}}+ b \sqrt{\frac{1-t}{2}}r\cos\theta
\end{equation}
and
\begin{equation}\label{mon4b}
|w(v,y,a)|^2   = a^2 \frac{1+t}{2} + b^2\frac{1-t}{2}r^2  + ab\sqrt{1-t^2}r\cos\theta\ ,
\end{equation}
and we deduce from (\ref{atl2})  that
\begin{eqnarray}\label{atl5a}
K_{a,\ell}h(t) &=&
\left(\frac{1}{1+t}\right)^{\ell/2}\int_0^1\int_0^\pi  h\left(\left[ a^2(1+t) + b^2(1-t) r^2 + 2ab\sqrt{1-t^2} r\cos\theta    \right]-1\right)\nonumber\\
&\times&\left| a^2(1+t) + b^2(1-t) r^2 + 2ab\sqrt{1-t^2} r\cos\theta    \right|^{\ell/2}\nonumber\\
&\times &P_\ell^{(\beta)}
\left(\frac{a\sqrt{1+t}  + b\sqrt{1-t}r\cos\theta}
{\left|a^2(1+t) + b^2(1-t) r^2 + 2ab\sqrt{1-t^2} r\cos\theta    \right|^{1/2}}
\right)
{\rm d}m_{\alpha,\beta}(r,\theta)\ .\nonumber\\
\end{eqnarray}

The  Laplace formula for the ultraspherical polynomials \cite{Szego}, p. 94, which is a simple consequence of Gegenbauer's identity, can be written as
\begin{equation}\label{Laplace}
P^{(\beta)}_\ell(x) = \frac{\Gamma(\beta+1/2)}{\sqrt \pi \Gamma(\beta)} \int_0^\pi \left( x + \sqrt{x^2-1}\cos \phi\right)^\ell
\sin^{2\beta -1}(\phi) d\phi \ .
\end{equation}
With 
${\displaystyle
x =\frac{ w(v,y,a)\cdot \widehat e}{|w(v,y,a)|}\ ,
}$
we obtain, using the binomial formula,
\begin{eqnarray}
& &\left| a^2(1+t) + b^2(1-t) r^2 + 2ab\sqrt{1-t^2} r\cos\theta    \right|^{\ell/2}\nonumber\\
&\times &P_\ell^{(\beta)}
\left(\frac{a\sqrt{1+t}  + b\sqrt{1-t}r\cos\theta}
{\left|a^2(1+t) + b^2(1-t) r^2 + 2ab\sqrt{1-t^2} r\cos\theta    \right|^{1/2}}
\right) \\
&=& \sum_{j=0}^\ell \left( \begin{array}{c} \ell \\
j 
\end{array}\right) (a\sqrt{1+t})^{\ell-j}(br\sqrt{1-t})^j  \\
&\times&\frac{\Gamma(\beta+1/2)}{\sqrt \pi \Gamma(\beta)} \int_0^\pi
(\cos \theta + \sqrt{\cos^2 \theta -1} \cos \phi)^j \sin^{2\beta-1} \phi d \phi \\
&=& \sum_{j=0}^\ell \left( \begin{array}{c} \ell \\
j 
\end{array}\right) (a\sqrt{1+t})^{\ell-j}(br\sqrt{1-t})^j  P_j^{(\beta)}(\cos \theta)
\end{eqnarray}
\lanbox

\noindent
{\bf Proof of Theorem \ref{betterprod}:} By choosing $h= p^{\alpha, \beta+\ell}$ Gasper's product formula (Theorem \ref{betterprod}) follows immediately from Lemma \ref{better}.
\lanbox

Our next goal is to show that  for all $\alpha> \beta > -1/2$, $K_{a,\ell}$ possesses the crucial properties of self-adjointness, polynomial preservation and the limiting value identity that it inherits from ${\mathcal  K}_a$ when $\alpha$ and $\beta$ satisfy (\ref{ab}).

\begin{lm}\label{gaspa2} For all $a\in (-1,1)$,  $\alpha>\beta>-1/2$, and integers $\ell\ge 0$, the operator
$K_{a,\ell}$  on ${\mathcal  C}([-1,1])$ as defined in (\ref{s44}) has the following properties:
\smallskip

\noindent{\it (1)} $K_{a,\ell}$ is self adjoint on $L^2(\mu^{\alpha,\beta})$.

\smallskip

\noindent{\it (2)} The space of polynomials of any fixed degree is invariant under $K_{a,\ell}$.

\smallskip

\noindent{\it (3)} For any continuous function $h$,
${\displaystyle \lim_{t\to 1} K_{a,\ell} h(t ) = a^\ell h(2a^2-1)}$.
\end{lm}

\noindent{\bf Proof:} It is obvious from (\ref{atl5a}) that 
$\lim_{t\to 1}K_{a,\ell}h(t)  = a^\ell h(2a^2-1)$,  property {\it (3)} is taken care of. Next consider the polynomial preservation, property{\it (2)}.
It suffices to show that for each natural number $n$, if $h(t) = (t+1)^n$, then $K_{a,\ell}h(t)$ is a polynomial of order $n$.

For this choice of $h$,
\begin{align*}
&h\left(\left[ a^2(1+t) + b^2(1-t) r^2 + 2ab\sqrt{1-t^2} r\cos\theta    \right]-1\right) =\\
& \sum_{m=1}^n \frac{n!}{(n-m)!m!}\left(a^2(1+t)+ b^2(1-t)\right)^{n-m}
\left(2ab\sqrt{1-t^2}r\cos(\theta)\right)^m\ .
\end{align*}
Thus, from Lemma~\ref{better}, 
$K_{a,\ell}h(t)$ is a sum of multiples of terms of the form
$$ Q(t) (1-t^2)^{m/2}  \left(\frac{1-t}{1+t}\right)^{k/2}\int_0^1\int_0^\pi  r^{m+k} \cos^m\theta  P_k^{(\beta)}(\cos \theta) {\rm d} m_{\alpha,\beta}\ ,$$
where $Q(t)$ is a polynomial of degree $n-m$. 
Then, be the orthogonality properties of the ultraspherical polynomials, 
$$\int_0^\pi  \cos^m\theta  P_k^{(\beta)}(\cos \theta)\sin^{2\beta}\theta{\rm d}\theta = 0$$
unless $m+k$ is even and $m\ge k$,  in which case
${\displaystyle (1-t^2)^{m/2}  \left(\frac{1-t}{1+t}\right)^{k/2} = \frac{(1-t)^{(m+k)/2}}{(1+t)^{(m-k)/2}}
}$
is a polynomial of degree $m$. Thus, for this choice of $h$,  $K_{a,\ell}h(t)$ is a sum of terms each of which is a polynomial of degree $n$, and
thus {\it (2)} is proved.

We next deal with self-adjointness. To see this in a simple way, we do not use the formula for $K_{a,\ell}$ given in Lemma~\ref{better}, but instead work directly from the expression (\ref{atl2}).
We shall show  that the bilinear form
\begin{equation}\label{atl3}
q(h_1,h_2) := 2c_{\alpha,\beta} \int_0^1h_1(2s^2-1)\left(K_{a,\ell}h_2\right)(2s^2-1) (1-s^2)^\alpha s^{2\beta+2\ell+1}{\rm d}s
\end{equation}
is symmetric.  This is easily seen in case   $\alpha$ and $\beta$ related to $m$ and $N$ through
(\ref{ab}) since then with $f_j(v) = h_j(2|v|^2-1)H(v)$, $j=1,2$, easy computations reveal that
the right hand side is a constant multiple of
$\langle f_1, {\mathcal  K}_a f_2\rangle_{L^2(\nu_{m,N})}$
To see this in general, we proceed {\em exactly} as in the proof of Lemma~\ref{gaspa}, making the same sequences of coordinate changes
$$(r,\theta) \to  (x,y) \to (x',y') \to  (\rho,\phi)\ .$$
Under this sequence of changes of variables, $|w(v,y,a)|$ becomes simply $\rho$, as we have seen in the proof of 
Lemma~\ref{gaspa}, and $w(v,y,a)\cdot \widehat e/|w(v,y,a)|$ becomes simply $\cos(\phi)$, as simple computations reveal.
Then, with $q(h_1,h_2)$ defined in (\ref{atl3}), we find that 
\begin{eqnarray}\label{atl7}
q(h_1,h_2)  &=& 2c_{\alpha,\beta}b^{-2\alpha} \int^1_0\int^1_{0} \int^{\pi}_0
h_1(2s^2-1)h_2(2\rho^2-1)\nonumber\\
&\times& (b^2-s^2-\rho^2 +2a\rho s\cos{\phi})_+^{\alpha-\beta-1}P_{\ell}^{(\beta)}(\cos\phi) \rho^{2\beta + \ell+1}\,
s^{2\beta + \ell+1}\sin^{2\beta}\phi{\rm d} \phi\, {\rm  d}\rho {\rm d}s\ ,\nonumber\\
\end{eqnarray}
This takes care of property ({\it 1}). \lanbox

With this lemma in hand, we now easily extend Lemma \ref{higherell}

\begin{thm}\label{higherell2} For all $\alpha>\beta > -1/2$, all integers $\ell\ge 0$,  and each $a\in (-1,1)$, the operator $K_{a,\ell}$
is self adjoint on $L^2(\mu^{(\alpha,\beta+\ell)})$, and is diagonalized by the Jacobi polynomial basis
$\{p_n^{(\alpha,\beta+\ell)}\}_{n\ge 0}$. Moreover, the corresponding sequence of
eigenvalues $\{\lambda_{n,\ell}(a)\}_{n\ge 0}$ is given by 
\begin{equation}\label{lrat2}
\lambda_{n,\ell}(a)= a^{\ell}\frac{p_n^{(\alpha,\beta+\ell)}(2a^2-1)}{p_n^{(\alpha,\beta+\ell)}(1)}\ .
\end{equation}
\end{thm}

\noindent{\bf Proof:} This is an immediate consequence of Lemma~\ref{gaspa2} and Theorem~\ref{dec2}. \lanbox.

\noindent{\bf Proof of Theorem ~\ref{genko}:}  By Theorem~\ref{higherell2},
\begin{equation}\label{mon5}
 \frac{p_n^{(\alpha,\beta+\ell)}(t)}{p_n^{(\alpha,\beta+\ell)}(1)}
p_n^{(\alpha,\beta+\ell)}(2a^2-1)=
a^{-\ell}\left(K_{a,\ell}p_n^{(\alpha,\beta+\ell)}\right)(t)\ .
\end{equation}

Now, the left hand side is a polynomial in $a$, and so is the right hand side. Hence we may extend
the range of $a$ from $[-1,1]$ to all of $\R$. Since in (\ref{atl6}), $b$ stands for $\sqrt{1-a^2}$, 
All odd terms in $\sqrt{1-a^2}$ must drop out of when the integration is made, and  for $a>1$, we will get the signs right if
we replace $b = \sqrt{1-a^2}$ with $i\sqrt{a^2-1}$.

Doing this, and then dividing both sides of (\ref{mon5}) by $(2a^2)^n$, and taking the limit $a\to \infty$, only the leading terms in the Jacobi polynomials contribute,
and we obtain
we obtain
\begin{align*}
 \frac{p_n^{(\alpha,\beta+\ell)}(t)}{p_n^{(\alpha,\beta+\ell)}(1)}
=
\int_0^1\int_0^\pi  &\left[\frac{(1+t) - (1-t)r^2}{2} + i\sqrt{1-t^2}r\cos\theta \right]^{n}\\
\times &\left[\sum_{k = 0}^{\ell} 
 \left( \begin{array}{c} \ell \\
k 
\end{array}\right)\left(\frac{1-t}{1+t}\right)^{k/2} (ir)^k P_k^{(\beta)}(\cos \theta)\right]
{\rm d}m_{\alpha,\beta}(r,\theta)\ .
\end{align*}
This is the formula in Theorem~\ref{genko}. \lanbox


\section{Bounds on the extremal eigenvalues and convergence of\\ associated eigenfuncton expansions} \label{GSP6}

Our objective in this section is to obtain bounds on the magnitudes of the eigenvalues in the extremal Markov sequences that govern the way these magnitudes decrease to zero as $n$ increases. We start with the
case of the ultraspherical polynomials.

\begin{thm}\label{ultrabound} For all $\gamma>0$, and all $-1 < a < 1$, 
\begin{equation}\label{ev1}
\left| \frac{p_n^{(\gamma)}(a)}{p_n^{(\gamma)}(1)}\right|\le
\frac{2c_{\gamma-1/2}}{(1-a^2)^{\gamma}}\left(\frac{n}{2}\right)^{-\gamma}\ .
\end{equation}
Moreover,   if $p>1/\gamma$,  $(K_a)^p$, the $p$th power of $K_a$,  is trace class. This criterion for belonging to the trace class is sharp 
in that for $a=0$, where exact calculuations are simple, one finds
$K_0^p$ is trace class if and only if  $p>1/\gamma$.
\end{thm}

We shall prove an analog of this Theorem for Jacobi polynomials, and the proof will be quite similar.
Therefore, before plunging into the details, we explain the strategy.

The starting point is the Laplace identity (\ref{Laplace}) which can be written as
$$\frac{p_n^{(\gamma)}(a)}{p_n^{(\gamma)}(1)}  =
 \int_{-1}^1  \left(a   + is \sqrt{1-a^2}\right)^n{\rm d}\mu^{(\gamma-1/2)}(s)\ .
 $$
Observe  that
\begin{equation}\label{L2}
|a   + is \sqrt{1-a^2}|^2 = 1 - (1-a^2)(1-s^2) \le 1 \ .
\end{equation}

Fixing $a$, define $C_\lambda$ to be the subset of $[-1,1]$ on which 
$|a   + is \sqrt{1-a^2}|^2 \ge 1-\lambda$.  It follows from (\ref{Laplace}) and layer--cake that
\begin{equation}\label{L1}
\left| \frac{p_n^{(\gamma)}(a)}{p_n^{(\gamma)}(1)}\right| \le 
 \int_{-1}^1  \left|a   + is \sqrt{1-a^2}\right|^n{\rm d}\mu^{(\gamma-1/2)}(s) \le
\frac{n}{2} \int_0^1 (1-\lambda)^{(n-2)/2}\mu^{(\gamma-1/2)}(C_\lambda){\rm d}\lambda\ .
\end{equation}

Hence, an estimate on the rate that $\mu^{(\gamma-1/2)}(C_\lambda)$ decreases to zero as $\lambda$
decreases to zero yields a bound on the rate at which $|p_n^{(\gamma)}(a)/p_n^{(\gamma)}(1)|$ decreases as $n$ increases.  This will yield us bounds that hold uniformly in $a$ in any compact subset of $(-1,1)$.
While we are ignoring phase cancelations in the estimate (\ref{L1}),  there are no phase cancelations for $a=0$, and an exact calculation gives the same $n^{-\gamma}$ decay. Thus no better bound can hold uniformly in $a$ on closed symmetric intervals of $(-1,1)$.

We prove a bound on $\mu^{(\gamma-1/2)}(C_\lambda)$ in the next lemma, and then proceed with the
proof of the theorem.

\begin{lm}\label{ultralem} \quad 
${\displaystyle \mu^{(\gamma-1/2)}(C_\lambda) \le 2c_{\gamma-1/2}\left(\frac{\lambda}{1-t^2}\right)^{\gamma}}\ .$
\end{lm}

\noindent{\bf Proof:} Note that from (\ref{L2}),
$|t   + is \sqrt{1-t^2}|^2 \ge 1-\lambda \iff  (1-t^2)(1-s^2) \le \lambda$. 
Hence, for $s\in  C_\lambda$, $1-s^2 < \lambda/(1-t^2)$, and therefore,
$$\mu^{(\gamma-1/2)}(C_\lambda) = 2c_{\gamma-1/2}\int_{\sqrt{1 - \lambda/(1-t^2)}}^1
(1-s^2)^{\gamma-1}{\rm d}s \le 2c_{\gamma-1/2}\left(\frac{\lambda}{1-t^2}\right)^{\gamma}.$$ \lanbox

\noindent{\bf Proof of Theorem~\ref{ultrabound}:}
Applying  Lemma~\ref{ultralem}  in (\ref{L1}), we  obtain
\begin{equation}\label{L3}
\left| \frac{p_n^{(\gamma)}(t)}{p_n^{(\gamma)}(1)}\right| \le
nc_{\gamma-1/2}\int_0^1 (1-\lambda)^{(n-2)/2}\left(\frac{\lambda}{1-t^2}\right)^{\gamma}{\rm d}\lambda\ = \frac{c_{\gamma-1/2}}{(1-t^2)^{\gamma}} n\frac{\Gamma(\frac{n}{2})\Gamma(\gamma+1)}{\Gamma(\frac{n}{2}+\gamma+1)}.
\end{equation}
Then since $(1-e^{-s})^{\gamma}\le s^{\gamma}$ for $\gamma\ge0$,
\begin{equation}\label{betaest}
 \frac{\Gamma(\frac{n}{2})}{\Gamma(\frac{n}{2}+\gamma+1)}=\frac{1}{\Gamma(\gamma+1)}\int_0^{\infty} e^{-\frac{n}{2}s}(1-e^{-s})^{\gamma}ds\le \left(\frac{n}{2}\right)^{-(\gamma+1)} \ .
\end{equation}
Combining this with (\ref{L3}) we obtain the bound (\ref{ev1}).

Finally, we consider the case $a=0$. Then there is no phase cancelation, and one readily computes
$$
\left|\frac{p_{2n}^{(\gamma)}(0)}{p_{2n}^{(\gamma)}(1)}\right|=
\frac{c_{\gamma-1/2}}{\sqrt{2}}\frac{\Gamma(n+\frac{1}{2})\Gamma(\gamma)}{\Gamma(n+\gamma+\frac{1}{2})}\sim n^{-\gamma}\ ,
$$
and we see that the $n^{-\gamma}$ bound on the absolute value of the extremal 
eigenvalues is the best possible  that can hold uniformly for $a$ in closed intervals of
$(-1,1)$. \lanbox

We now turn to the analog of Theoem~\ref{ultrabound} for Jacobi polynomials.

\begin{thm}\label{eigs} For all $\alpha> \beta > -1/2$, all $a\in (-1,1)$,  all  $n>0$, and all $\ell\ge 0$,
\begin{equation}\label{con11}
\left| \frac{p_n^{(\alpha,\beta+\ell)}(a)}{p_n^{(\alpha,\beta+\ell)}(1)}\right| \le 
\left[1 +  \left(\frac{1-a}{1+a}\right)^{1/2}  \right]^\ell K_{\alpha,\beta}(a) 
\Gamma\left(\alpha+\frac{3}{2}\right)   \left(\frac{n}{2}\right)^{-(\alpha+1/2)} \ ,
\end{equation}
where 
\begin{equation}\label{con}
K_{\alpha, \beta}(a) = c_{\alpha,\beta}\pi^{-1}2^{2(\alpha - \beta+2)}(1-a)^{-(2\alpha-\beta + 1)}
(1+a)^{-(\beta+1/2)}\ .
\end{equation}
Thus, for $p>1/(\alpha+1/2)>0$, $(K_{a,\ell})^p$ is trace class.
\end{thm}

\begin{remark}  The exponent on $n$ is determined by $\alpha$ alone; it is independent of
$\beta$ and $\ell$.
\end{remark}

We begin with a lemma  that is the analog of Lemma~\ref{ultralem}:

\begin{lm}\label{size} For each fixed $a \in (-1,1)$ and $\lambda> 0$, 
define $C_\lambda$ to be the
subset of $[0,1]\times [0,\pi]$ given by
\begin{equation}\label{clamdef}
C_\lambda := \{(r,\theta)\  : \  \ R^2(r,\theta) > 1- \lambda\ \}\ .
\end{equation}
where
$$R(r,\theta) =  \left|\frac{(1+a) - (1-a)r^2}{2} + i\sqrt{1-a^2}r\cos\theta \right|\ .$$
Then 
the measure of $C_\lambda$ with respect to $\d  m_{\alpha,\beta}$ satisfies
\begin{equation}\label{cset}
m_{\alpha,\beta}(C_\lambda) \le
K_{\alpha, \beta}(a)\lambda^{\alpha+1/2}\ .
\end{equation}
where $K_{\alpha, \beta}(a)$ is given by (\ref{con}).
\end{lm}

\noindent{\bf Proof:} 
Define $A = (1+a)/2$ and $B = (1-a)/2$.  Then we can write
$$R^2(r,\theta) = A^2 + B^2r^2 + 2ABr^2\cos(2\theta)\ .$$
Since $A+B =1$,
$1 - R^2(r,\theta) =  B^2(1-r^4) + 2AB(1 - r^2\cos(2\theta))$.
This is a sum of positive terms, and so  for any $\lambda>0$, whenever 
$1 - R^2(r,\theta) \le \lambda$, we have both 
\begin{equation}\label{twob}
 B^2(1-r^4) < \lambda \qquad{\rm and}\qquad 2AB(1 - r^2\cos(2\theta)) < \lambda\ .
 \end{equation}
The first of these conditions implies
$1 - r^2 < \lambda/(B^2(1+r^2)(<  \lambda/B^2$, and then $1-r < (1-r)(1+r) = 1-r^2 < \lambda/B^2$, so 
that $r > 1 -\lambda/B^2$. Thus,everywhere on $C_\lambda$,
$$1-r^2 < \frac{\lambda}{B^2}\qquad{\rm and}\qquad  r > 1 -\frac{\lambda}{B^2}\ .$$

Next, we turn to the second condition in (\ref{twob}). This can be written as
$r^2\cos(2\theta) > 1 - \lambda/(2AB)$, which certainly implies
$\cos(2\theta) > 1 - \lambda/(2AB)$, which implies that
$\sin^2\theta < \lambda/(4AB)$. Finally, since on $[0,\pi/2]$, $(2/\pi)\theta \le \sin(\theta)$, with
a similar estimate on $[\pi/2,\pi]$, the second condition in (\ref{twob}) 
implies that 
either $0 \le \theta \le  (\pi/4)\sqrt{\lambda/(AB)}$, or else $\pi -  (\pi/4)\sqrt{\lambda/(AB)} \le \theta \le \pi$. 

Altogether then, 
$\{(r,\theta)\  : \ R^2(r,\theta) > 1- \lambda\ \}$
is contained in 
$[1 - \lambda/B^2,1]\times [0,(\pi/4)\sqrt{\lambda/(AB)}] \ \cup  [1 - \lambda/B^2,1]\times
 [\pi -(\pi/4)\sqrt{\lambda/(AB)}, \pi] $,
and moreover, everywhere on this set,
$$1-r^2 \le  \frac{\lambda}{B^2} \qquad{\rm and}\qquad \sin^2\theta \le 
\frac{\lambda}{4AB}\ .$$
Integrating over the two rectangles using the above bounds yields the estimate
$$m_{\alpha,\beta}(C_\lambda) \le c_{\alpha,\beta}\pi^{-1}2^{3-2\beta}  B^{-(2\alpha-\beta+1/2)}
A^{-(\beta +1/2)}\lambda^{\alpha+1/2}\ .$$
Replacing $A$ and $B$ by their definition in terms of $t$, one obtains the bound (\ref{cset}).
 \lanbox

\noindent{\bf Proof of Theorem~\ref{eigs}:} The starting point is Theorem~\ref{genko},
which provides the identity
\begin{align*}
 \frac{p_n^{(\alpha,\beta+\ell)}(a)}{p_n^{(\alpha,\beta+\ell)}(1)}
=
\int_0^1\int_0^\pi  &\left[\frac{(1+a) - (1-a)r^2}{2} + i\sqrt{1-a^2}r\cos\theta \right]^{n}\\
\times &\left[\sum_{k = 0}^{\ell} 
 \left( \begin{array}{c} \ell \\
k 
\end{array}\right)\left(\frac{1-a}{1+a}\right)^{k/2} (ir)^k P_k^{(\beta)}(\cos \theta)\right]
{\rm d}m_{\alpha,\beta}(r,\theta)\ .
\end{align*}
By the definition of $R(r,\theta)$ in Lemma~\ref{size}, and the fact that $P_n^{(\gamma)}(x) \le 1$,
we obtain
\begin{eqnarray}\label{jacrat}
\left| \frac{p_n^{(\alpha,\beta+\ell)}(t)}{p_n^{(\alpha,\beta+\ell)}(1)}\right|
&\le& 
\int_0^1\int_0^\pi  R(r,\theta)^n
\left[\sum_{k = 0}^{\ell} 
 \left( \begin{array}{c} \ell \\
k 
\end{array}\right)\left(\frac{1-a}{1+a}\right)^{k/2} r^k \right]
{\rm d}m_{\alpha,\beta}(r,\theta)\nonumber\\
&=& 
\int_0^1\int_0^\pi  R(r,\theta)^n
\left[1 + 
\left(\frac{1-a}{1+a}\right)^{1/2} r \right]^\ell
{\rm d}m_{\alpha,\beta}(r,\theta)\nonumber\\
&\le&\left[1 +  \left(\frac{1-a}{1+a}\right)^{1/2}  \right]^\ell
\int_0^1\int_0^\pi  R(r,\theta)^n
{\rm d}m_{\alpha,\beta}(r,\theta)\nonumber\\
&=&\left[1 +  \left(\frac{1-a}{1+a}\right)^{1/2}  \right]^\ell
\frac{n}{2} \int_0^1 (1-\lambda)^{(n-2)/2}m_{\alpha,\beta}(C_\lambda){\rm d}\lambda\ .\nonumber\\
\end{eqnarray}

Now applying Lemma ~\ref{size}, and then estimating the ratio of Gamma functions as in (\ref{betaest}),
\begin{eqnarray}
\left| \frac{p_n^{(\alpha,\beta+\ell)}(a)}{p_n^{(\alpha,\beta+\ell)}(1)}\right| &\le& 
\left[1 +  \left(\frac{1-a}{1+a}\right)^{1/2}  \right]^\ell\frac{n}{2}K_{\alpha,\beta}(a)  \int_0^1  (1- \lambda)^{(n-2)/2} 
\lambda^{\alpha+1/2}{\rm d}\lambda\nonumber\\
&=&\left[1 +  \left(\frac{1-a}{1+a}\right)^{1/2}  \right]^\ell\frac{n}{2}K_{\alpha,\beta}(a)  \frac{\Gamma(\frac{n}{2})\Gamma(\alpha+\frac{3}{2})}{\Gamma(\frac{n}{2}+\alpha+\frac{3}{2})}\nonumber\\
&\le&\left[1 +  \left(\frac{1-a}{1+a}\right)^{1/2}  \right]^\ell K_{\alpha,\beta}(a) 
\Gamma\left(\alpha+\frac{3}{2}\right)   \left(\frac{n}{2}\right)^{-(\alpha+1/2)}  \nonumber\\
\end{eqnarray}
\lanbox

We may apply these results to study the convergence of the eigenfunction expansions for the
operators $K_a$ and $K_{a,\ell}$.  Indeed, since eigenvalues of $K_{a,\ell}$ are
$$\lambda_{n,\ell}(a)= a^{\ell}\frac{p_n^{(\alpha,\beta+\ell)}(2a^2-1)}{p_n^{(\alpha,\beta+\ell)}(1)}\ ,$$
and the eigenfunctions are the $p_n^{(\alpha,\beta+\ell)}$, the formal eigenfunction expansion of
the Kernel for  $K_{a,\ell}$ is
$$a^\ell \sum_{n=0}^\infty \frac{p_n^{(\alpha,\beta+\ell)}(2a^2-1) p_n^{(\alpha,\beta+\ell)}(x)p_n^{(\alpha,\beta+\ell)}(y)}
{p_n^{(\alpha,\beta+\ell)}(1)}\ .$$
The eigenvalue bounds obtained above can be used to show that for $a\in (-1,1)$ and $\alpha>1/2$, these formal series actually
converge uniformly for $x$ and $y$ in compact intervals of $(-1,1)$.  To do this, we need bounds on the
eigenfunctions as well as the eigenvalues. Because of the close relation between the eigenfunctions and the eigenvalues in this context, we could obtain the eigenfunction bouts from Theorem~\ref{eigs},
but it will be instructive to obtain these instead from a well known but deep result of
Nevai, Erdelyi, and Magnus. \cite{NEM}:    For all
 $\alpha\ge -1/2$ and  $\beta\ge
  -1/2$ and all non negative integers $n$, 
\begin{equation}\label{nembnd}
{\rm max}_{x\in[-1,1]}\sqrt{1-x^2}w(x)p^{\alpha,\beta}_n(x)^2\le\frac{2e(2+\sqrt{\alpha^2+\beta^2})}{\pi},
\end{equation}
Thus, for each $a\in (-1,1)$ and each $r<1$, there is a constant $C$ such that
$$\frac{p_n^{(\alpha,\beta)}(2a^2-1) p_n^{(\alpha,\beta)}(x)p_n^{(\alpha,\beta)}(y)}
{p_n^{(\alpha,\beta)}(1)} \le Cn^{-(\alpha+1/2)}$$
uniformly for $x,y\in [-r,r]$.  With the $\alpha = \beta$ cases coming from Theorem~\ref{ultrabound} and
(\ref{convent}), this proves:

\begin{thm} For all $\alpha> 1/2$ and $\alpha\ge  \beta > -1/2$, and all $-1 <x,y,z < 1$, the sum
\begin{equation}\label{gassum}
\sum_{n=0}^\infty \frac{p_n^{(\alpha,\beta+\ell)}(x)p_n^{(\alpha,\beta+\ell)}(y)p_n^{(\alpha,\beta+\ell)}(z)}
{p_n^{(\alpha,\beta+\ell)}(1)}
\end{equation}
converges absolutely and uniformly on compacts, and the operator whose kernel the sum defines is
trace class. 
\end{thm}

For $\ell = 0$, this is the eigenfunction expansion of Gasper's operator $K_{a,0}$, which is Markov.
Thus, for  $\alpha> 1/2$ and $\alpha\ge  \beta > -1/2$, where the sum in (\ref{gassum}), 
converges pointwise, it defines a kernel that is pointwise positive..

\section{Historical Remarks}\label{h}

In this section we give a brief discussion of work done on the Markov sequence problem for Jacobi polynomials by Gasper, Koornwinder and Askey with the aim of clarifying the context of the present paper.

When Gasper took up his work on the Markov sequence problem for Jacobi polynomials, the main obstacle was the lack of an analog to Gegenbauer's indentity. Therefore, Gasper worked backwards towards one:  It is clear that the  operator $K_{a,0}$ on $L^2(\mu^{(\alpha,\beta)})$ defined by
\begin{equation}\label{gasper2a}
K_{a,0}\psi(x) = \sum_{n=0}^\infty \frac{p_n^{(\alpha,\beta)}(2a^2-1)}
{p_n^{(\alpha,\beta)}(1)} p_n^{(\alpha,\beta)}(x) \int_{-1}^1p_n^{(\alpha,\beta)}(y)\psi(y){\rm d}\mu^{(\alpha,\beta)}(y)
\end{equation}
for $a\in (-1,1)$ is self adjoint, satisifes $K_{a,0}1=1$, and has $\{p_n^{(\alpha,\beta)}(z)/p_n^{(\alpha,\beta)}(1)\}_{n\ge 0}$ as its sequence of eigenvalues, so that if there is to be an analog of Gegenbauer's identity
for Jacobi polynomials, it {\em must} refer to this operator. 

What is not at all clear from the eigenfunction expansion is whether or not $K_{a,0}$ preserves positivity, or whether $K_{a,0}$ even has a kernel
$K_{a,0}(x,y)$, which formally would be
\begin{equation}\label{gasper2c}
K_{a,0}(x,y) =  \sum_{n=0}^\infty \frac{p_n^{(\alpha,\beta)}(2a^2-1) p_n^{(\alpha,\beta)}(x) p_n^{(\alpha,\beta)}(y)}{p_n^{(\alpha,\beta)}(1)}\ .
\end{equation}

Gasper's Theorem as stated above was proved by him in  \cite{G1}, where he
 evaluated the sum. This is rather involved, but here is a brief sketch:
In \cite{G1} Gasper formally defines  a kernel $G(x,y,z;\alpha,\beta)$ as
$$
G(\cos2\phi,\cos2\psi,\cos2\theta;\alpha,\beta)=
\sum_{n=0}^{\infty}h_n^{\alpha,\beta}\frac{p_n^{\alpha,\beta}(\cos2\phi)}{p_n^{\alpha,\beta}(1)}p_n^{\alpha,\beta}(\cos2\psi) p_n^{\alpha,\beta}(\cos2\theta),
$$
where $h_n^{\alpha,\beta}$ is the square of the inverse of the $L^2$ norm 
of $p_n^{\alpha,\beta}$. Then using a formula of Watson, \cite{W} p. 413,  he shows that
this  sum of triple products of Jacobi polynomials is equal to an integral of a triple product of Bessel functions with the restriction $\alpha>-1/2,  \beta>-1/2$ and $\cos\theta\ne|\cos(\psi\pm\phi)|$. When $\alpha>\beta>-1/2$, he was able to evaluate  the  integral of the triple product of Bessel functions 
with the result that
\begin{align}\label{gone}
&G(\cos2\phi,\cos2\psi,\cos2\theta;\alpha,\beta)\nonumber\\&=\frac{\Gamma(\alpha+1)(\sin\phi\sin\psi\sin\theta)^{-2\alpha}}{2^{\alpha+\beta+1}\Gamma(\alpha-\beta)\Gamma(\beta+1/2)\Gamma(1/2)}\int_0^A(1-\cos^2\phi-\cos^2\psi-\cos^2\theta+\nonumber\\&2\cos\phi\cos\psi\cos\theta\cos\gamma)^{\alpha-\beta-1}\sin^{2\beta}\gamma d\gamma
\end{align}
where $A$ is $0$, ${\displaystyle {\rm arccos}\left(\frac{\cos^2\phi+\cos^2\psi+\cos^2\theta-1}{2\cos\phi\cos\psi\cos\theta}\right)}$, or  $\pi$ depending on whether $\sin^2\phi\sin^2\psi$ is less than between, or greater than the two numbers $(\cos\phi\cos\psi\pm\cos\theta)^2$.
from this Gasper concludes that $G$ is  non-negative.
Then in \cite{G2}, using the evaluation of the triple integral in terms of hypergeometric functions, Gasper is able to show that
$G$ is nonnegative if $\alpha\ge\beta>-1,\alpha>-1/2$, and either
$\beta\ge-1/2$ or $\alpha+\beta\ge0$. Comparing this with equation
(\ref{prfo8}) we see that the kernel is positive so that Gasper's result lays the
foundation for a convolution structure associated with  Jacobi polynomials. 

Later Koornwinder \cite{korn} gave another proof of 
Gasper's Theorem in the case $\alpha>\beta>-1/2$. Here Koornwinder defines the 
kernel $G$ as the integral (\ref{s44}) given above, then he uses his Laplace type integal 
representation for Jacobi polynomials and duality to show that the kernel is 
equal to the triple sum of Jacobi polynomials. Koornwinder obtained his 
Laplace type formula using group theoretic methods and Askey \cite{A1} gave a 
simple analytic proof using Bateman's integral relation between hypergeometric 
functions. The fact that the kernel is continuous and of bounded variation allows Koornwinder to show, 
using the Dirichlet-Jordan test \cite{Zy} p. 57 and the equiconvergence of  
Jacobi series and cosine series \cite{Szego} p. 246, that for $\alpha>\beta>-1/2$ the triple sum 
converges uniformly on compact subsets of
$0<\phi,\psi,\theta<\frac{\pi}{2}$. Later Koornwinder and Schwartz
\cite{KooSchwa} extended these results polynomials orthogonal on the biangle, 
triangle, and simplex. 

In contrast our approach is much more  in the spirit  of  Bakry and Mazet \cite{BM}, in which 
they solved the Markov sequence problem for various systems of orthogonal polynomials by applying functional analytic techniques to 
certain well chosen kernels of self--adjoint operators. 
Likewise, our approach starts with the construction of the family of operators 
the operators ${\mathcal  K}_a$.  The motivation for considering the family ${\mathcal  K}_a$ comes from previous work on the Kac model \cite{CCL2}, \cite{CGL}. In particular, for the restricted parameter
values discussed in Section~\ref{Bl}, ${\mathcal  K}_a$ is easily seen to be self adjoint, to preserve
polynomials and to enjoy the evaluation property, and hence, by Theorem~\ref{dec2}  its eigenvalues
can be expressed as ratios of Jacobi polynomials. Restricting this operator to the radial functions
gives us the operator $K_{a,0}$, at least for the half integral values of $\alpha$ and $\beta$
given in (\ref{ab}). In
a further step we extend the operators  $K_{a,0}$ to the full range, and obtain Gasper's kernel.
$K_{a,0}$ is shown to be an extremal Markov operator from which Gasper's product formula follows

Moreover, the same can be done with the operators
$K_{a,\ell}$ that appear as restrictions of
the operators ${\mathcal  K}_a$ to the other invariant angular momentum subspaces, and in this way we obtain Theorems ~\ref{betterprod} and \ref{genko}.

One final remark is that the convolution structure for Jacobi polynomials 
can be used to show that Jacobi polynomials form a strong 
polynomial hypergroup (see Bloom and Heyer \cite{BlHe} for a definition of 
hypergroups). A sequence of orthonormal polynomials $\{f_n\}$ with $f_0=1$ 
such that \ref{product} holds is said to have the hypergroup property 
(see Bakry and Huet \cite{BH}). A result of Connett and Schwartz
\cite{CoS} (see also Theorem 3.4 in \cite{BlHe}) essentially shows that the 
only unit orthogonal polynomials systems in one dimension that satisfy 
\ref{product} are, 
up to a translation and scaling,  the Jacobi polynomials. This can be
proved  by showing that 
if a sequence of orthonormal polynomials satisfy \ref{product} then after 
translation and scaling they also satisfy the differential equation
satisfied by the Jacobi polynomials.  The argument does not, however, provide another proof of
Gasper's Thoerem asserting that 
the Jacobi
polynomials  do indeed satisfy a product formula.

\appendix


\section{The parabolic biangle and triangle polynomial product formulas}\label{parabiangle}
The {\em parabolic biangle}
is the region
\begin{equation}\nonumber
B=\{(x_1,x_2): 0 \le x_2^2 \le x_1 \le 1\} \ ,
\end{equation}
and the polynomials $r_{n,m}^{\alpha, \beta}(x_1,x_2)$ that are orthogonal on $B$ with respect to the measure
\begin{equation}\nonumber
(1-x_1)^\alpha(x_1-x_2^2)^{\beta-1/2} dx_1 dx_2 
\end{equation}
can be written, with our conventions as
\begin{equation}\label{bianglepoly}
r_{n,m}^{\alpha, \beta}(x_1,x_2) = x_1^{\frac{n}{2}} p^{(\alpha, \beta+n)}_{m}(2x_1-1) p_n^{(\beta)}(\frac{x_2}{\sqrt x_1}) \ .
\end{equation}
The total degree is $n+m$.

The product formula Koornwinder and Schwartz \cite{KooSchwa} is then given by
\begin{equation}\label{biangleprod}
\frac{r_{n,m}^{\alpha, \beta}(x_1^2,x_2) r_{n,m}^{\alpha, \beta}(y_1^2,y_2)}{ r_{n,m}^{\alpha, \beta}(1,1)}
=\int_{I,J^3} r_{n,m}^{\alpha, \beta}(E^2,EG) d \nu^{\alpha,\beta}(r,t_1,t_2,t_3)
\end{equation}
where
$$
d\nu^{\alpha,\beta}(r,t_1,t_2,t_3) = d\mu^{\beta}(t_2)d\mu^{\beta}(t_3)dm_{\alpha,\beta}(r,t_1) \ ,
$$
where $\mu^{\beta}$ is given by \ref{mualpha}, $m_{\alpha,\beta}$ by
\ref{dmalpha} with $t_1=\cos{\theta}$.  
Here $I= [0,1], J = [-1,1]$.
The symbols $E$ and $G$ are given by
$$
E= \left(x_1^2 y_1^2 +(1-x_1^2)(1-y_1^2)r^2 +2x_1y_1(1-x_1^2)^{1/2}(1-y_1^2)^{1/2}r t_1 \right)^{1/2} \ ,
$$
$$
G = D(C,D(x_2/x_1,y_2/y_1;1, t_2); 1,t_3)
$$
where
$$
C=\frac{D(x_1,y_1;r,t_1)}{E(x_1,y_1;r,t_1)}
$$
and, generally,
$$
D(x,y;r,t) = xy +(1-x^2)^{1/2}(1-y^2)^{1/2}r t \ .
$$
In order to prove this product formula (\ref{biangleprod}) using the method developed here,
we must analyze the operator 
\begin{equation}
({\mathcal K}_{y_1,y_2}h)(x_1,x_2) := \int_{I,J^3} h(E^2,EG) d \nu^{\alpha,\beta}(r,t_1,t_2,t_3) 
\end{equation}
on the Hilbert space ${\mathcal H}$ given by the inner product
$(\cdot,\cdot)$ defined above. Direct calculations that lead to the proof
of self-adjointness of the operator ${\mathcal K}_{y_1,y_2}$ defined above
seem to be very involved. A substantial simplification is achieved by
writing ${\mathcal K}_{y_1,y_2}$ in terms of Gegenbauer operators. For the
triangle case the Gasper operator will be used instead.
It is an immediate consequence of the following theorem.
\begin{thm}
The operator  ${\mathcal K}_{y_1,y_2}$ is a selfadjoint linear operator on ${\mathcal H}$. It is positivity
preserving, preserves the function $1$ and the space of polynomials in the two variables $x_1,x_2$
of a given degree. Further we have the evaluation formula
\begin{equation}
\lim_{(x_1,x_2) \to (1,1), (x_1,x_2) \in B} ({\mathcal K}_{y_1,y_2}h)(x_1,x_2) = h(1,1) \ .
\end{equation}
\end{thm}

\noindent{\bf Proof:}
Clearly, ${\mathcal K}_{y_1,y_2}$ is positivity preserving and preserves the function $1$. The evaluation formula
follows by noting that $E$ and $G$ tend to $1$ as $(x_1, x_2)$ tend to $(1,1)$ in the biangle. To see the statement 
concerning polynomial preservation it suffices to prove it for a general monomial $x_1^m x_2^n$. That is we have to show that
\begin{equation}\nonumber
\int_{I,J^3} E^{2m}(EG)^n  d \nu^{\alpha,\beta}(r,t_1,t_2,t_3)
\end{equation}
is a polynomial in the variables $(x_1,x_2)$ of total degree less than or equal $2m+n$. We shall use the fact that the measure $d \mu^{\alpha, \beta}$ is
even in $t_1,t_2,t_3$.
Now, 
\begin{eqnarray}
G^n &=&\left[CD(x_2/x_1,y_2/y_1;1, t_2) + (1-C^2)^{1/2}(1-D(x_2/x_1,y_2/y_1;1, t_2)^2)^{1/2} t_3\right]^n  \nonumber \\ 
&=& \sum_{k=0}^n  \left( \begin{array}{c} n \\ k  \end{array}\right) C^{n-k} D(x_2/x_1,y_2/y_1;1, t_2)^{n-k}
(1-C^2)^{k/2}(1-D(x_2/x_1,y_2/y_1;1, t_2)^2)^{k/2} t_3^k \nonumber
\end{eqnarray}
and integrating this expression with respect to the $t_3$ we see that only terms with even $k$ contribute and we obtain an expression
of the form
\begin{equation}\nonumber
 \sum_{k=0}^{[n/2]} c(n,2k) C^{n-2k} D(x_2/x_1,y_2/y_1;1, t_2)^{n-2k}
(1-C^2)^{k}(1-D(x_2/x_1,y_2/y_1;1, t_2)^2)^{k}
\end{equation}
where $c(n,2k)$ are positive coefficients. This expression can be rewritten as
\begin{equation}
 \sum_{k=0}^{[n/2]} c(n,2k) C^{n-2k} (1-C^2)^{k}  \sum_{p=0}^k \left( \begin{array}{c} k \\ p \end{array}\right)
 (-1)^p D(x_2/x_1,y_2/y_1;1, t_2)^{n-2k +2p} \ .
\end{equation}
Applying the binomial formula to the expression $D(x_2/x_1,y_2/y_1;1, t_2)^{n-2k +2p}$ and integrating with respect to $t_2$
leaves us with a polynomial in the variables $(x_2/x_1)^2$ if $n$ is even or is
of the form $x_2/x_1$ times a polynomial in $(x_2/x_1)^2$ otherwise. Moreover, it has degree not larger than $n$.
The remaining terms are, when multiplied by $E^{2m+n}$, of the form
\begin{eqnarray}
& &E^{2m+n}C^{n-2k} (1-C^2)^{k} = E^{2m} D(x_1,y_1;r,t_1)^{n-2k}(E^2 -D(x_1,y_1;r,t_1)^2)^k \nonumber \\
&=& E^{2m} D(x_1,y_1;r,t_1)^{n-2k}((1-x_1^2)(1-y_1^2) r^2 (1-t_1^2))^k
\end{eqnarray}
which when integrated over $t_1$ yields a polynomial in $x_1^2$ if $n$ is even
or is of the form $x_1$ times a polynomial on $x_1^2$ otherwise. It has degree not larger than $2m+n$.
Thus, after performing the integration over the variables $r,t_1,t_2,t_3$ one obtains a sum of terms of the form
$$
(x_1^2)^q \left(\frac{x_2^2}{x_1^2}\right)^r = x_1^{2q-2r}x_2^{2r} \ n \ {\rm even} \ , \ {\rm with}  \ 2q \le 2m+n
$$
or
$$
x_1(x_1^2)^q \frac{x_2}{x_1}\left(\frac{x_2^2}{x_1^2}\right)^r = x_1^{2q-2r}x_2^{2r+1} \ n \ {\rm even} \ , \ {\rm with}  \ 2q+1 \le 2m+n \ .
$$
Thus we obtain a polynomial of the form $p(x_1^2,x_2)$ whose total degree is not larger than the one we started with.
It remains to show selfadjointness. It is convenient to write the inner product
\begin{equation}
(f,g) := \int_0^1 dx_1(1-x_1)^\alpha  \int_{-\sqrt{x_1}}^{\sqrt{x_1}}dx_2 (x_1-x_2^2)^{\beta-1/2} f(x_1,x_2) g(x_1,x_2) \ .
\end{equation}
in terms of the functions
$$
F(\rho, s) := f(\rho^2, \rho s) \ ,  \  G(\rho, s) := g(\rho^2, \rho s) \ ,
$$
\begin{equation}
(f,g) = 2 \int_0^1 d \rho \rho^{2 \beta +1} (1-\rho^2)^\alpha \int_{-1}^1 ds (1-s^2)^{\beta-1/2}
F(\rho, s) G(\rho,s)  =: \langle F,G\rangle \ .
\end{equation}
This follows from the definition of $(\cdot,\cdot)$ by a simple change of variables. Thus instead of the variables
$x_1,x_2$ we have the variables $x_1 = \rho^2$ and $x_2 = \rho s$. Likewise we write
$y_1 = y^2$ and $y_2 =y t$. Note, that in this notation the form of the orthogonal polynomials (\ref{bianglepoly})
becomes apparent.

The operator ${\mathcal K}_{y_1,y_2}$ in these variables is given by
\begin{equation}\label{op}
[{\mathcal K}_{y,t} H](\rho , s) = \int_{I,J^3} d \nu^{\alpha,\beta}(r,t_1,t_2,t_3) H(E, G)
\end{equation}
where
\begin{equation}\nonumber
E = E(\rho, y; r, t_1) = (\rho^2 y^2 +(1-\rho^2)(1-y^2) r^2 + 2\rho y (1-\rho^2)^{1/2}
(1-y^2)^{1/2} r t_1 )^{1/2}
\end{equation}
and
\begin{equation}\nonumber
G = G(\rho, s, y, t: r, t_1, t_2, t_3)
= D(C, D(s,t;1,t_2); 1, t_3) \ ,
\end{equation}
\begin{equation}\nonumber
C = \frac{D(\rho , y;r, t_1)}{  E(\rho, y; r, t_1)} =:\frac{D}{E} \ .
\end{equation}
As before, 
\begin{equation}\nonumber
D(a,b;r,t) = ab +(1-a^2)^{1/2}(1-b^2)^{1/2} r t \ .
\end{equation}
Recall that
$$
K_a f(t) = \int_{-1}^1 f(at+s \sqrt{1-a^2}\sqrt{1-t^2})d\mu^{\beta-1/2}(s)
$$
which was used for the Gegenbauer product formula. Now note that
\begin{equation}\nonumber
[K^{(2)}_tK^{(2)}_{\frac{D}{E}}H](E,s) =  \int_{-1}^1
d\mu^{\beta-1/2}(t_2) \int_{-1}^1 d\mu^{\beta-1/2} (t_3) H(E, G) \ .
\end{equation}
The superscript  $(2)$ indicates that the operator acts on the second variable of the function.
Now
\begin{eqnarray}
& &\langle F, [{\mathcal K}_{y,t} H] \rangle =  2 \int_0^1 \int_{-1}^1 d m_{\alpha, \beta}(r, t_1)
 \int_0^1 d \rho \rho^{2 \beta +1} (1-\rho^2)^\alpha \nonumber\\
& & \times \int_{-1}^1 ds (1-s^2)^{\beta-1/2}
F(\rho, s) [K^{(2)}_tK^{(2)}_{\frac{D}{E}}H](E,s)   \ .\nonumber
\end{eqnarray}
Now we proceed using Koornwinder's change of variables: First by going to cartesian
coordinates we get
\begin{eqnarray}
& &\int_0^1 \int_{-1}^1 d m_{\alpha, \beta}(r, t_1) 
[K^{(2)}_tK^{(2)}_{\frac{D}{E}}H](E,s)  \nonumber\\
& =& l_{\alpha,\beta}  \int_{-\infty}^\infty du \int_0^\infty dv  (1-u^2 -v^2)_+^{\alpha-\beta-1}
v^{2\beta} [K^{(2)}_tK^{(2)}_{\frac{D}{E}}H](E,s)\nonumber
\end{eqnarray}
where $E$ and $D$ expressed in these new variables are given by
\begin{equation}\nonumber
E = \left[(1-\rho^2)(1-y^2) v^2 +[(1-\rho^2)^{1/2}(1-y^2)^{1/2}u + \rho y]\right]^{1/2} \ , \ 
D= (1-\rho^2)^{1/2}(1-y^2)^{1/2}u + \rho y  \ .
\end{equation}
Here $l_{\alpha,\beta}$ is the normalizing constant.
By scaling the variables $(u,v) \to (1-\rho^2)(1-y^2)]^{1/2} (u,v)$we get
\begin{eqnarray}
& &\int_0^1 \int_{-1}^1 d m_{\alpha, \beta}(r, t_1) 
[K^{(2)}_tK^{(2)}_{\frac{D}{E}}H](E,s)  \nonumber\\
& =& l_{\alpha,\beta} [(1-\rho^2)(1-y^2)]^{-\alpha}  \int_{-\infty}^\infty du \int_0^\infty dv  ((1-\rho^2)(1-y^2)-u^2 -v^2)_+^{\alpha-\beta-1} \nonumber\\ 
&\times &v^{2\beta}[K^{(2)}_tK^{(2)}_{\frac{D}{E}}H](E,s)\nonumber
\end{eqnarray}
Shifting $u \to u+\rho y$ and reverting to polar coordinates we obtain
\begin{eqnarray}
& &\int_0^1 \int_{-1}^1 d m_{\alpha, \beta}(r, t_1) 
[K^{(2)}_tK^{(2)}_{\frac{D}{E}}H](E,s) \nonumber \\
& =& l_{\alpha,\beta} [(1-\rho^2)(1-y^2)]^{-\alpha}  \int_0^\infty dr r^{2\beta + 1}  \int_{-1}^1 d\sigma  (1-\sigma^2)^{\beta-1/2} (1-\rho^2-y^2-r^2 + 2\rho y r \sigma )_+^{\alpha-\beta-1}\nonumber \\ 
&\times &[K^{(2)}_tK^{(2)}_{\frac{D}{E}}H](E,s)\nonumber
\end{eqnarray}
where this time
\begin{equation}\nonumber
E = r \ , \ D = r\sigma \ .
\end{equation}
Collecting the terms we obtain
\begin{eqnarray}
& &\langle F, [{\mathcal K}_{y,t} H] \rangle\nonumber  \\
&=& 2 l_{\alpha,\beta} (1-y^2)^{-\alpha}  \int_0^1 d \rho \rho^{2 \beta +1}  \nonumber \\
& \times & \int_0^\infty dr r^{2\beta + 1}  \int_{-1}^1 d\sigma  (1-\sigma^2)^{\beta-1/2} (1-\rho^2-y^2-r^2 + 2\rho y r \sigma )_+^{\alpha-\beta-1} \nonumber\\ 
&\times & \int_{-1}^1 ds (1-s^2)^{\beta-1/2}  F(\rho, s) [K^{(2)}_tK^{(2)}_\sigma H](r,s)\nonumber
\end{eqnarray}

The operators $K^{(2)}_t$ and $K^{(2)}_\sigma$ are selfadjoint with respect to
the scalar product with the measure $ds (1-s^2)^\beta$. Moreover,
they commute which follows from the fact that they have the same
eigenvectors the ultraspherical polynomials. 
Hence
\begin{eqnarray}
& &\langle F, [{\mathcal K}_{y,t} H] \rangle \nonumber \\
&=& 2 c_{\alpha,\beta} (1-y^2)^{-\alpha}  \int_0^1 d \rho \rho^{2 \beta +1}  \nonumber \\
& \times & \int_0^\infty dr r^{2\beta + 1}  \int_{-1}^1 d\sigma  (1-\sigma^2)^{\beta-1/2} (1-\rho^2-y^2-r^2 + 2\rho y r \sigma )_+^{\alpha-\beta-1} \nonumber\\ 
&\times & \int_{-1}^1 ds (1-s^2)^{\beta-1/2}  [K^{(2)}_tK^{(2)}_\sigma F](\rho, s) H(r,s)\nonumber
\end{eqnarray}
which, since the expression is symmetric in $\rho$ and $r$, equals
$\langle [{\mathcal K}_{y,t} F], H \rangle$.
\lanbox

An similar argument gives the product formula for triangle polynomials
first derived by Koornwinder and Schwartz. Recall the scalar product for the orthogonal polynomials on the triangle:

\begin{equation}
(f,g) := \int_0^1 dx_1(1-x_1)^\alpha  \int_{0}^{x_1} (x_1-x_2)^\beta x_2^{\gamma}f(x_1,x_2) g(x_1,x_2) \ .
\end{equation}

The polynomials, orthogonal in this inner product, are denoted by
$$
R^{\alpha, \beta,\gamma}_{n,k}(x_1,x_2)=R_{n-k}^{\alpha,\beta+k+1/2}(2x_1-1)x_1^k R_k^{\beta,\gamma}(2\frac{x_2}{x_1}-1).
$$
The product formula is
\begin{equation}
R^{\alpha, \beta,\gamma}_{n,k}(x^2_1,x^2_2)R^{\alpha, \beta,\gamma}_{n,k}(y^2_1,y^2_2)=\int_{I^4\times J^3}R^{\alpha, \beta,\gamma}_{n,k}(E^2,E^2 H^2)d\mu^{\alpha,\beta,\gamma}
\end{equation}
where
 $$d\mu^{\alpha,\beta,\gamma}(r_1,r_2,r_3,r_4,\psi_1,\psi_2,\psi_3)\nonumber\\
=dm_{\beta,\gamma}(r_4,\psi_3)dm_{\beta,\gamma}(r_3,\psi_4)d\nu^{\beta,\gamma-1/2}(r_2)dm_{\alpha,\beta+\gamma+1}(r_1,\psi_1)\ ,$$
with $dm_{\alpha<\beta}$  given by equation~(\ref{dmalpha}), 
\begin{equation}\label{nualpha}
d\nu^{\beta, \gamma-1/2}(r_2)=\hat c_{\beta,\gamma}(1-r_2)^{\beta}r_2^{\gamma-1/2}dr_2.
\end{equation}
and where
\begin{eqnarray}
 H &=& H(x_1,x_2,x_3,x_4;r_1,r_2,r-2,r_3,\psi_1,\psi_2,\psi_3)\nonumber\\ &=&E\left([1-r_2)C^2+r_2]^{1/2},E(\frac{x_2}{x_1},\frac{y_2}{y_1};r_3,\psi_2);r_4,\psi_3\right)\ .\nonumber
\end{eqnarray}
The definitions for $E$, $D$ and $C$ are as in the biangle formula.

It is convenient to rewrite the inner product in terms of the functions
$$
F(x_1, x_2) := f(x_1^2, x_1^2 x_2^2) \ ,  \  G(x_1, x_2) := g(x_1^2, x_1^2 x_2^2) \ ,
$$
so the inner product is given by
$$
(f,g) = 4 \int_0^1 dx_1 x_1^{2 (\beta+\gamma) +3} (1-x_1^2)^\alpha \int_{0}^1 dx_2 x_2^{2\gamma+1}(1-x_2^2)^\beta
F(x_1, x_2) G(x_1,x_2)  =: \langle F,G\rangle \ .
$$
This follows from the definition of $(\cdot,\cdot)$ by a simple change of variables.
We also let $(y_1, y_2)\to(y_1^2,y_1^2 y_2^2)$.
The triangle product formula of Koornwinder and Schwartz is now given in terms of the following
operator:
\begin{eqnarray}\label{opt}
& &[{\mathcal G}_{y_1,y_2} G](x_1 , x_2) = \int_0^1 \int_{-1}^1 d m_{\alpha, \beta+\gamma+1/2}(r_1,t_1)\int_{0}^1 d\nu^{\beta,\gamma-1/2}(r_2)\nonumber\\ & &  \int_{0}^1\int_{-1}^1
dm_{\beta,\gamma} (r_3,t_2) \int_{0}^1\int_{-1}^1 dm_{\beta,\gamma} (r_4,t_3) G(E, H)\nonumber
\end{eqnarray}
where
$$
E = E(x_1, y_1; r_1, t_1) = (x_1^2 y_1^2 +(1-x_1^2)(1-y_1^2) r_1^2 + 2x_1 y_1 (1-x_1^2)^{1/2}
(1-y_1^2)^{1/2} r_1 t_1 )^{1/2}
$$
$$
C = \frac{D(x_1 , y_1;r_1, t_1)}{  E(x_1, y_1; r_1, t_1)},
$$
where generally
$$
D(a,b;r,t) =ab  +(1-a^2)^{1/2}(1-b^2)^{1/2} r t \ .
$$

Gasper's operator can be rewritten as
$$
K^{\alpha,\beta}_{y,0} f(t) = \int_0^1\int_{-1}^1 f((y^2t^2+(1-t^2)^{1/2}r_2^2+2yt(1-t^2)^{1/2}(1-y^2)^{1/2}r_2 t_3))^{1/2} dm_{\alpha,\beta}(r_2,t_3),
$$
Thus
$$
[K^{\beta,\gamma,(2)}_{y_2,0}K^{\beta,\gamma,(2)}_{[(1-r_2)C^2+r_2]^{1/2},0}G](E,x_2) =  \int_0^1\int_{-1}^1
dm_{\beta,\gamma}(r_2,t_3)\int_0^1 \int_{-1}^1 dm_{\beta,\gamma} (r_3,t_4) G(E, H) \ .
$$
As above the superscript  $(2)$ indicates that the operator acts on the second variable of the function. In the formuls below the constant $c_{\alpha,\beta,\gamma}$
denotes the products of the various constants normalizing the measures that we use.

Thus the inner product can be written as,
\begin{eqnarray}
& &\langle F, [{\mathcal G}_{y_1,y_2} G] \rangle =  4\int_0^1(1-x_1^2)^{\alpha} x_1^{2(\beta+\gamma)+3}\int_0^1 d\nu^{\beta,\gamma-1/2}(r_2) \int_0^1 \int_{-1}^1 d m_{\alpha, \beta+\gamma+1/2}(r_1, t_1) \nonumber\\
& & F(x_1, x_2) [K^{\beta,\gamma,(2)}_{y_2,0}K^{\beta,\gamma,(2)}_{[(1-r_2)C^2+r_2]^{1/2},0}G](E,x_2)   \ .\nonumber
\end{eqnarray}

Now we proceed using Koornwinders change of variables and following the discussion of
the biangle formula we see that $E$ and $D$ become
$E = r \ , \ D = r\sigma$. Collecting terms and making the final change of varialbes $u=[(1-r_2)\sigma^2+r_2]^{1/2}$ we obtain

\begin{eqnarray}
& &\langle F, [{\mathcal G}_{y_1,y_2} G] \rangle \nonumber \\
&=& 8 c_{\alpha,\beta,\gamma} (1-y^2)^{-\alpha}  \int_0^1 dx_1 x_1^{2(\beta+\gamma) +3} \int_0^\infty dr r^{2(\beta+\gamma) + 3}\int_{-1}^1 d\sigma (1-\sigma^2)^{\beta+\gamma+1/2}\times \nonumber \\&
&\int_{\sigma}^1 du(1-u^2)^{\beta}(u^2-\sigma^2)^{\gamma-1/2}u (1-x_1^2-y_1^2-r^2 + 2x_1 y_1 r \sigma )_+^{\alpha-\beta-\gamma-2} \nonumber \\ 
&\times & \int_0^1 dx_2(1-x^2_2)^{\beta}x_2^{2\gamma+1} F(x_1, x_2) [K^{\beta,\gamma,(2)}_{y_2,0}K^{\beta,\gamma, (2)}_{u,0} G](r,x_2)\nonumber 
\end{eqnarray}
The proof of Gasper's theorem shows that the operators $K^{\beta,\gamma,(2)}_{y_2,0}$ and $K^{\beta,\gamma, (2)}_{u,0}$ are selfadjoint with respect to
the scalar product with respect to the measure $dx_2 (1-x_2^2)^{\beta}x_2^{2\gamma+1}$. Moreover,
they commute which follows from the fact that they have the same eigenvectors, i.e.,
the Jacobi Polynomials. 
Hence the self adjointness follows as in the biangle formula. 

The polynomial preservation also follows from an argument similar to the
biangle formula. Here we get even powers of $E$ and $H$ and use the fact that the integrals
over $t_1...t_4$ are symmetric.

\end{document}